\documentclass[journal]{IEEEtran}
\usepackage{etex}
\usepackage[vlined,ruled]{algorithm2e}
\usepackage[dvipsnames]{xcolor}
\usepackage{pgfpages}
\usepackage[cmex10]{amsmath}		
\usepackage{amssymb, mathrsfs}
\usepackage{cite}
\usepackage{url}

\usepackage{dsfont}
\usepackage{siunitx}
\usepackage{verbatim}									
\usepackage{mathtools}									
\usepackage{siunitx}									
\usepackage{mdframed}

\usepackage{graphicx}
\usepackage{epstopdf}

\usepackage[font=small]{caption}
\usepackage{subcaption}
\usepackage[activate={true,nocompatibility},final,tracking=true,kerning=true,spacing=true,factor=1100,stretch=10,shrink=10]{microtype}

\usepackage{array}
\usepackage{multirow}
\usepackage{enumerate}
\usepackage{booktabs}

\graphicspath{{./fig/}}

\usepackage{pgfplots}
\pgfplotsset{width=8cm,compat=newest}
\newcommand*{\damping}{0.006}%
\newcommand*{\freq}{25}%
\pgfmathsetmacro{\freqd}{sqrt(1-(\damping)^2)*\freq}%

\pgfplotsset{
    standard/.style={
    axis x line=middle,
    axis y line=middle,
    enlarge x limits=0.15,
	enlarge y limits=0.15,
	every axis plot post/.style={mark options={fill=black}},
	}
}

\pgfplotsset{%
    ,compat=1.12
    ,every axis x label/.style={at={(current axis.right of origin)},anchor=north west}
    ,every axis y label/.style={at={(current axis.above origin)},anchor=north east}
    }

\usepackage{tikz}
\usetikzlibrary{arrows}
\usetikzlibrary{decorations.pathmorphing}
\usetikzlibrary{decorations.markings}
\usetikzlibrary{snakes}
\usetikzlibrary{matrix}
\usetikzlibrary{calc}
\usetikzlibrary{patterns}
\usetikzlibrary{positioning}
\usetikzlibrary{shapes}
\usetikzlibrary{shapes.symbols}
\usetikzlibrary{shapes.geometric}
\usetikzlibrary{shadows.blur}
\usetikzlibrary{shapes.arrows}
\usetikzlibrary{shapes.multipart}
\usetikzlibrary{shapes.callouts}
\usetikzlibrary{shapes.misc}

\tikzstyle{every node}=[font=\small]
\tikzstyle{every path}=[line width=0.8pt,line cap=round,line join=round]

\usepackage[american,smartlabels]{circuitikz}
\ctikzset{bipoles/thickness=1}
\ctikzset{bipoles/length=0.8cm}
\ctikzset{bipoles/diode/height=.375}
\ctikzset{bipoles/diode/width=.3}
\ctikzset{tripoles/thyristor/height=.8}
\ctikzset{tripoles/thyristor/width=1}
\ctikzset{bipoles/vsourceam/height/.initial=.7}
\ctikzset{bipoles/vsourceam/width/.initial=.7}

\newcommand{\real}{\mathbb{R}}

\newcommand{\until}[1]{\{1,\dots, #1\}}

\DeclareMathOperator*{\maximize}{maximize} 									
\DeclareMathOperator*{\argmin}{argmin}

\newcommand{\vect}[1]{\mathbbold{#1}}
\newcommand{\vones}[1][]{\vect{1}_{#1}}
\newcommand{\vzeros}[1][]{\vect{0}_{#1}}
\DeclareSymbolFont{bbold}{U}{bbold}{m}{n}
\DeclareSymbolFontAlphabet{\mathbbold}{bbold}

\newcommand{\map}[3]{#1: #2 \rightarrow #3}

\newcommand\oprocendsymbol{\hbox{$\square$}}
\newcommand\oprocend{\relax\ifmmode\else\unskip\hfill\fi\oprocendsymbol}

\newtheorem{theorem}{Theorem}[section]
\newtheorem{lemma}[theorem]{Lemma}

\newtheorem{corollary}[theorem]{Corollary}
\newtheorem{proposition}[theorem]{Proposition}

\newenvironment{pfof}[1]{\vspace{1ex}\noindent{\itshape Proof of
    #1:}\hspace{0.5em}} {\hfill\oprocend\vspace{1ex}}

\newcommand{\diag}{\mathrm{diag}}

\title{Input-Output Performance of Linear-Quadratic Saddle-Point Algorithms with Application to Distributed Resource Allocation Problems%
\thanks{
}}

\author{John W. Simpson-Porco, Bala Kameshwar Poolla, Nima Monshizadeh, and Florian D\"{o}rfler%
    \thanks{%
    J. W. Simpson-Porco is with the Department of Electrical and Computer Engineering, University of Waterloo, ON, Canada. Email: {\tt jwsimpson@uwaterloo.ca}. N. Monshizadeh is with the Engineering and Technology Institute, University of Groningen, The Netherlands. Email: {\tt n.monshizadeh@rug.nl}.
     B. K. Poolla and F. D{\"o}rfler are with the Automatic Control Laboratory, Swiss Federal Institute of Technology (ETH) Z\"urich, Switzerland. Email: {\tt \{bpoolla,dorfler\}@ethz.ch}. The work of J. W. Simpson-Porco is supported in part by the NSERC Discovery Grant RGPIN-2017-04008. The work of B. K. Poolla and F. D\"orfler is supported by ETH start-up funds and the SNF AP Energy Grant \#160573. %
}}

\begin{document}
\maketitle
\thispagestyle{empty}
\pagestyle{empty}

\begin{abstract}
Saddle-point or primal-dual methods have recently attracted renewed interest as a systematic technique to design distributed algorithms which solve convex optimization problems. When implemented online for streaming data or as dynamic feedback controllers, these algorithms become subject to disturbances and noise; convergence rates provide incomplete performance information, and quantifying input-output performance becomes more important. We analyze the input-output performance of the continuous-time saddle-point method applied to linearly constrained quadratic programs, providing explicit expressions for the saddle-point $\mathcal{H}_2$ norm under a relevant input-output configuration. We then proceed to derive analogous results for regularized and augmented versions of the saddle-point algorithm. We observe some rather peculiar effects -- a modest amount of regularization significantly improves the transient performance, while augmentation does not necessarily offer improvement. We then propose a distributed dual version of the algorithm which overcomes some of the performance limitations imposed by augmentation. Finally, we apply our results to a resource allocation problem to compare the input-output performance of various centralized and distributed saddle-point implementations and show that distributed algorithms may perform as well as their centralized counterparts.
\end{abstract}


\section{Introduction}
\label{Section: Introduction}

Saddle-point methods are a class of continuous-time gradient-based algorithms for solving constrained convex optimization problems.
Introduced in the early 1950s \cite{TK:56,KA-LH-HU:58}, these algorithms are designed to seek the saddle points of the optimization problem's Lagrangian function. 
These saddle points are in one-to-one correspondence with the solutions of the first-order optimality (KKT) conditions, and the algorithm therefore drives its internal state 
towards the global optimizer of the convex program; see \cite{DF-FP:10,JW-NE:11,AC-EM-JC:15} for convergence results. 

Recently, these algorithms have attracted renewed attention for e.g., in the context of machine learning \cite{KCK-SM-RP:18}, in the control literature for solving \emph{distributed} convex optimization problems \cite{BG-JC:14}, where agents cooperate through a communication network to solve an optimization problem with minimal or no centralized coordination. Applications of distributed optimization include utility maximization \cite{DF-FP:10}, congestion management in communication networks \cite{SHL-DEL:99}, and control in power systems \cite{FD-JWSP-FB:13y,NL-CL-ZC-SHL:13,AC-JC:14,EDA-SVD-GBG:15,TS-CDP-AvdS:17,EM-CZ-SL:17}.
While most standard optimization algorithms require centralized information to compute the optimizer, saddle-point algorithms often yield distributed strategies in which agents perform state updates using only locally measured information and communication with some subset of other agents. We refer the reader to \cite{JTW-MA:04, JW-NE:10,JW-NE:11, GH-HD-MBE:14, GD-ME:14, TH-NC-TI-NL:16, HDN-TLV-KT-JS:18} for control-theoretic interpretations of these algorithms.

Rather than solve the optimization problem offline, it is desirable to run these distributed algorithms online as controllers, in feedback with system and/or disturbance measurements, to provide references so that the optimizer can be tracked in real-time as operating conditions change. Such algorithms offer promise for online optimization, especially in scenarios with streaming data.
However, when saddle-point methods are implemented online as controllers, they become subject to disturbances arising from fluctuating parameters and noise (the precise nature of these disturbances being application dependent). The standard method for assessing optimization algorithms \textemdash{} namely, convergence rate analysis \textemdash{} is now insufficient to capture the performance of the algorithm. Indeed, an algorithm with a fast convergence rate would be inappropriate for control applications if it responded poorly to disturbances during transients, or if it greatly amplifies measurement noise in steady-state. 

The appropriate tool for measuring dynamic algorithm performance is instead the system norm, as commonly used in feedback system analysis to capture system response to exogenous disturbances. Recent work in this direction includes input-to-state-stability results \cite{AC-EM-SL-JC:16,AC-EM-SL-JC:17}, finite $\mathscr{L}_2$-gain analysis \cite{JWSP:16f}, and the robust control framework proposed in \cite{LL-BR-AP:16,BH-LL:17}. The purpose of this paper is to continue this line of investigation. In particular, the case of saddle-point algorithms, applied to optimization problems with quadratic objective functions and linear equality constraints, leads to a very tractable instance of this analysis problem where many basic questions can be asked and accurately answered. Relevant questions include
\begin{enumerate}
\item[(i)] how do saddle-point algorithms amplify disturbances which may enter the objective function and/or equality constraints?
\item[(ii)] how does performance in the presence of disturbances change when the initial optimization problem is reformulated (e.g., dual or distributed formulations)?
\item[(iii)] how do standard modifications to optimization algorithms, such as regularization and Lagrangian augmentation, affect these results?
\end{enumerate}
\smallskip

\emph{Contributions:} The three main contributions of this paper are as follows.\footnote{A preliminary version of these results with application to power system control appeared in the conference article \cite{JWSP-BKP-NM-FD:16c}. In contrast to the conference article, this paper reports the proofs of Theorems \ref{Thm:PrimalOutputs} and \ref{Thm:AugmentedPrimalOutputs}, studies the effect of regularization (Theorem \ref{Thm:PrimalOutputs_reg}), studies both centralized and distributed dual saddle-point approaches (Proposition \ref{Prop:Dual-Ascent} and Corollary \ref{Cor:ADD}),  and studies the application of these results to resource allocation problems (Section \ref{Sec:Resource}).}
First, in Section \ref{Section: Performance} we consider the effect of disturbances on the saddle-point dynamics arising from linearly constrained, convex quadratic optimization problems. We quantify the input-output performance of the method via the $\mathcal{H}_2$ system norm, and \textemdash{} for a relevant input-output configuration \textemdash{} derive an explicit expression for the norm as a function of the algorithm parameters. We find (Theorem \ref{Thm:PrimalOutputs}) that the squared $\mathcal{H}_2$ norm scales linearly with the number of disturbances to both the primal and dual variable dynamics. 

Second, we study two common modifications to the Lagrangian optimization paradigm: regularization and augmentation. We show that regularization strictly improves the transient $\mathcal{H}_2$ performance of saddle-point algorithms (Theorem \ref{Thm:PrimalOutputs_reg}). However, this improvement in performance is not usually monotone in the regularization parameter, and the system norm may achieve its global minimum at some finite regularization parameter value. For augmented Lagrangian saddle-point methods, we derive an explicit expression for the $\mathcal{H}_2$ norm (Theorem \ref{Thm:AugmentedPrimalOutputs}); the results show that augmentation may either improve or deteriorate the $\mathcal{H}_2$ performance. For cases when standard augmentation deteriorates performance, we propose augmented dual distributed saddle-point algorithm which strictly improves performance (Section \ref{Sec:DistributedDual}). 

Third and finally, in Section \ref{Sec:Resource} we apply our results to resource allocation problems, comparing and contrasting the input-output performance of the different algorithms we have considered. The results show that distributed implementations can perform equally well as centralized implementations, but that significant performance differences can appear between the algorithms once augmentation is considered. 

Taken together, these results provide fairly complete answers to the questions (i)--(iii) outlined in the introduction, for the class of problems considered. A similar study in a $\mathcal{H}_{\infty}$ or $\mathscr{L}_{2}$-gain framework requires a significantly different analysis, and is outside the scope of the present paper; see Section \ref{Section: Conclusions} for open directions. As background, in Section \ref{Sec:PrimalDualMethods} we review saddle-point algorithms for the relevant class of optimization problems, then recall the basic facts about the $\mathcal{H}_2$ norm as a measure of input-output system performance.

\smallskip

\emph{Notation:}  The $n \times n$ identity matrix is $I_n$, $\vzeros[]$ is a matrix of zeros of appropriate dimension, while $\vones[n]$ (resp. $\vzeros[n]$) are $n$-vectors of all ones (resp. zeros). If $\map{f}{\real^n}{\real}$ is differentiable, then $\map{\frac{\partial f}{\partial x}}{\real^n}{\real^n}$ is its gradient. For $A \in \real^{n \times n}$, $A^{\sf T}$ is its transpose and $\mathrm{Tr}(A) = \sum_{i=1}^n A_{ii}$ is its trace. If $S \in \real^{r \times n}$ has full row-rank, then $SS^{\dagger} = I_r$ where $S^{\dagger} = S^{\sf T}(SS^{\sf T})^{-1}$ is the Moore-Penrose pseudoinverse of $S$. For a positive semidefinite matrix $Q \succeq \vzeros[]$, $Q^{\frac{1}{2}}$ is its square root. The symbol $\otimes$ denotes the Kronecker product. Given elements $\{a_i\}_{i=1}^{n}$ (scalars, vectors, or matrices), $\mathrm{col}(a_1,\ldots,a_n) = (a_1^{\sf T},\ldots,a_{n}^{\sf T})^{\sf T}$ denotes the vertically concatenation of the elements (assuming compatible dimensions), and $\text{blkdiag}(a_1, a_2, \ldots, a_n)$ is a block matrix with the elements $\{a_i\}$ on the diagonals.

\emph{Graphs and graph matrices:} A graph is a pair $\mathcal{G} = (\mathcal{V},\mathcal{E}_{\rm u})$, where $\mathcal{V}$ is the set of vertices (nodes) and $\mathcal{E}_{\rm u}$ is the set of undirected edges (unordered pairs of nodes). The set of neighbours of node $i \in \mathcal{V}$ are denoted by $\mathcal{N}(i)$. 
If a label $e \in \until{|\mathcal{E}_{\rm u}|}$ and an arbitrary orientation is assigned to each edge, we can define a corresponding directed edge set $\mathcal{E} \subset \mathcal{V} \times \mathcal{V}$ with elements $e \sim (i,j) \in \mathcal{E}$. The \emph{node-edge incidence matrix} $E \in \real^{|\mathcal{V}|\times|\mathcal{E}|}$ is defined component-wise as $E_{ke} = 1$ if node $k$ is the source node of edge $e$ and as $E_{ke} = -1$ if node $k$ is the sink node of edge $e$, with all other elements being zero. If the graph is connected, then $\mathrm{ker}(E^{\sf T})=\mathrm{Im}(\vones[|\mathcal{V}|])$. A graph is a \emph{tree} (or \emph{acyclic}) if it contains no cycles, and in this case $\mathrm{ker}(E) = \{\vzeros[|\mathcal{E}|]\}$.

\section{Saddle-Point Methods and $\mathcal{H}_2$ Performance}
\label{Sec:PrimalDualMethods}

\subsection{Review of Saddle-Point Method}
\label{Sec:ReviewPrimal}

We consider the constrained quadratic optimization problem
\begin{equation}\label{Eq:ConvexEqualityConstrained}
\begin{aligned}
& \underset{x \in \mathbb{R}^{n_x}}{\text{minimize}}
& & J(x) := \frac{1}{2}x^{\sf T}Qx + x^{\sf T} c \\
& \text{subject to}
&& Sx = W_bb\,,
\end{aligned}
\end{equation}
where $x \in \real^{n_x}$, $Q = Q^{\sf T} \succ \vzeros[]$ is positive definite, $c \in \real^{n_x}$ and $b \in \real^{n_b}$ are parameter vectors, and $S \in \real^{n_r \times n_x}, W_{b} \in \real^{n_r \times n_b}$ with $n_r < n_x$. We make the blanket assumption that $S$ and $W_{b}$ have full row rank, which simply means that the constraints $Sx = W_b b$ are not redundant. While the right-hand side $W_{b} b$ of the constraints is apparently over-parameterized, this formulation is natural when considering particular problem instances. In the resource allocation problem of  Section \ref{Sec:Resource}, $W_{b} = \begin{bmatrix}
1 & 1 & \cdots & 1
\end{bmatrix}$ and $b$ is a vector of demands; the product $W_{b}b$ is simply the total demand.

The problem \eqref{Eq:ConvexEqualityConstrained} describes only a subclass of the optimization problems to which saddle-point algorithms are applicable; more generally one considers strictly convex costs and convex inequality constraints as well. We restrict our attention to \eqref{Eq:ConvexEqualityConstrained}, as this case will allow LTI system analysis techniques to be applied, and represents a large enough class of problems to yield some general insights. Intuitively, the performance of the saddle-point algorithm on \eqref{Eq:ConvexEqualityConstrained} should indicate a ``best'' case performance for the general case, as the objective $J(x)$ is smooth and strongly convex, and \eqref{Eq:ConvexEqualityConstrained} is free of hard inequality constraints. See Section \ref{Section: Conclusions} for further discussion.
Under these assumptions, the convex problem \eqref{Eq:ConvexEqualityConstrained} has a finite optimum, the equality constraints are strictly feasible, and \eqref{Eq:ConvexEqualityConstrained} may be equivalently studied through its Lagrange dual with zero duality gap \cite{SB-LV:04}. 
The Lagrangian $\map{L}{\real^{n_x} \times \real^{n_r}}{\real}$ of the problem \eqref{Eq:ConvexEqualityConstrained} is
\begin{equation}\label{Eq:Lagrangian}
L(x,\nu) = \frac{1}{2}x^{\sf T}Qx + c^{\sf T}x + \nu^{\sf T}(Sx-W_b b)\,,
\end{equation}
where $\nu \in \real^{n_r}$ is a vector of Lagrange multipliers. By strong duality, the KKT conditions
\begin{equation}\label{Eq:KKT}
\begin{aligned}
\frac{\partial L}{\partial x}(x,\nu) &= \vzeros[n_x] \quad \Longleftrightarrow \quad \vzeros[n_x] =  Qx + S^{\sf T}\nu + c\,,\\
\frac{\partial L}{\partial \nu}(x,\nu) &= \vzeros[n_r] \quad \Longleftrightarrow \quad \vzeros[n_r] = Sx - W_bb\,,
\end{aligned}
\end{equation}
are necessary and sufficient for optimality. From these linear equations one can quickly  compute the unique global optimizer $(x^{\star},\nu^{\star})$ to be

\begin{equation}\label{Eq:Equilibrium}
\begin{bmatrix}
x^{\star} \\ \nu^{\star}
\end{bmatrix} = 
\begin{bmatrix}
-Q^{-1}(S^{\sf T}\nu^{\star} + c )\\
-(SQ^{-1}S^{\sf T})^{-1}(W_b b+SQ^{-1}c)
\end{bmatrix}\,.
\end{equation}
While \eqref{Eq:Equilibrium} is the exact solution to the optimization problem \eqref{Eq:ConvexEqualityConstrained}, its evaluation requires centralized knowledge of the matrices $S, Q, W_b$ and the vectors $c$ and $b$. If any of these parameters change or evolve over time, the optimizer should be recomputed. In many multi-agent system applications, the cost matrix $Q$ is diagonal or block-diagonal and $J(x) = \sum_{i}\frac{q_i}{2}x_i^2+c_ix_i$ is therefore a sum of local costs. Finally, the constraints encoded in $S$ are often sparse, mirroring the topology of an interaction or communication network between agents. These factors motivate the solution of \eqref{Eq:ConvexEqualityConstrained} in an online distributed fashion, where agents in the network communicate and cooperate to calculate the global optimizer.

A simple continuous-time algorithm to seek the optimizer is the \emph{saddle-point} or \emph{primal-dual} method \cite{MB-GHG-JL:05,JW-NE:11,AN-AO:09b,DF-FP:10,AC-EM-JC:15}
\begin{align*}
{\mathcal{T}_x}\dot{x} &= -\frac{\partial}{\partial x} L(x,\nu)\,, \quad
{\mathcal{T}_\nu}\dot{\nu} = +\frac{\partial}{\partial \nu} L(x,\nu)\,,
\end{align*}
which here reduces to the affine dynamical system

\begin{subequations}\label{Eq:SaddlePointContinuous}
\begin{align}
{\mathcal{T}_x}\dot{x} &= -Qx - S^{\sf T}\nu - c\\
{\mathcal{T}_\nu} \dot{\nu} &= Sx-W_b b\,,
\end{align}
\end{subequations}

where ${\mathcal{T}_x, \mathcal{T}_\nu} \succ \vzeros[]$ are positive definite diagonal matrices of time-constants. By construction, the equilibrium points of \eqref{Eq:SaddlePointContinuous} are in one-to-one correspondence with the solutions of the KKT conditions \eqref{Eq:KKT}, and the system is internally exponentially stable \cite{AC-EM-JC:15}.

\smallskip

\begin{lemma}[\bf Global Convergence to Optimizer]\label{Lem:Stable}
The unique equilibrium point $(x^{\star},\nu^{\star})$ given in \eqref{Eq:Equilibrium} of the saddle-point dynamics \eqref{Eq:SaddlePointContinuous} is globally exponentially stable, with exponential convergence rate $\propto 1/\tau_{\rm max}$ where $\tau_{\rm max} = \max(\max_{i \in \until{n}}{\mathcal{T}_{x,ii}},\max_{i\in \until{r}} {\mathcal{T}_{\nu,ii}})$.
\end{lemma}

\smallskip
\begin{pfof}{Lemma \ref{Lem:Stable}}

The proof of Lemma follows by using the Lyapunov candidate $V(x,\nu) =0.5 (x-x^\star)^{\sf T}{\mathcal{T}_{x}}(x-x^{\star}) + 0.5(\nu-\nu^{\star})^{\sf T}{\mathcal{T}_{\nu}}(\nu-\nu^{\star})+ \varepsilon (\nu-\nu^{\star})^{\sf T} S {\mathcal{T}_{x}}(x-x^\star)$ for $\varepsilon>0$, i.e., consider the Lyapunov candidate
$$
P = \frac{1}{2}\begin{bmatrix}
\mathcal{T}_x & {\varepsilon} \mathcal{T}_xS^{\sf T}\\
{\varepsilon} S\mathcal{T}_x & \mathcal{T}_\nu
\end{bmatrix}\,,
$$
which is positive definite since for $\varepsilon$ sufficiently small. With $A$ as the system matrix in \eqref{Eq:SaddlePointContinuous}, we compute that
\begin{align*}
A^TP+PA &= -\begin{bmatrix}\Delta & \varepsilon QS^{\sf T}/2 \\ \varepsilon SQ/2 & \varepsilon SS^{\sf T} \end{bmatrix}\,,
\end{align*}
where $\Delta = Q-\frac{\varepsilon}{2}(\mathcal{T}_xS^{\sf T}\mathcal{T}_{\nu}^{-1}S+S^{\sf T}\mathcal{T}_{\nu}^{-1}S\mathcal{T}_x)$. Since $\varepsilon > 0$ and $S$ has full row-rank, $\varepsilon SS^{\sf T}$ is positive definite. Moreover, $\Delta$ is positive definite if $\varepsilon$ is sufficiently small. Standard Schur complement results then yield that $A^{\sf T}P + PA \prec \vzeros[]$ if and only if
$$
SS^{\sf T} - \frac{\varepsilon}{4} SQ\Delta^{-1}QS^{\sf T} > 0\,,
$$
which holds for $\varepsilon$ sufficiently small as $\lim_{\varepsilon \to 0}\Delta= Q$.

As the convergence rate is dictated by the eigenvalues, it follows from a simple variant of \cite[Theorem 3.6]{MB-GHG-JL:05} that the rate is $\propto$ $\min_{ii}(\text{blkdiag}({\mathcal{T}_x^{-1}, \mathcal{T}_\nu^{-1}}))$ or  $\propto 1/\max(\max_{i \in \until{n}}{\mathcal{T}_{x,ii}},\max_{i\in \until{r}} {\mathcal{T}_{\nu,ii}})$.
\end{pfof}

With stability settled, in what follows we will focus exclusively on quantifying transient input-output performance of \eqref{Eq:SaddlePointContinuous} in the presence of exogenous disturbances.

\subsection{System Performance in the $\mathcal{H}_2$ Norm}
Consider the linear time-invariant system
\begin{equation}\label{Eq:GeneralMIMO}
\begin{aligned}
\dot{x} &= Ax + B\eta\\
z &= Cx\,,
\end{aligned}
\end{equation}

where $\eta$ is the disturbance input signal and $z$ is the performance output. With $x(0) = \vzeros[]$, we denote the linear operator from $\eta$ to $z$ by $G$. If \eqref{Eq:GeneralMIMO} is input-output stable, its $\mathcal{H}_2$ norm $\|G\|_{\mathcal{H}_2}$ is defined as 
$$
\|G\|_{\mathcal{H}_2}^2 := \frac{1}{2\pi}\int_{-\infty}^{\infty} \mathrm{Tr}(G(-j\omega)^{\sf T}G(j\omega))\,\mathrm{d}\omega\,,
$$
where $G(j\omega) = C(j\omega I-A)^{-1}B$ is the frequency response of \eqref{Eq:GeneralMIMO}.

Another interpretation of $\|G\|_{\mathcal{H}_2}^2$ is as the steady-state variance of the output
\begin{equation}
\label{Eq:H2norm}
\|G\|_{\mathcal{H}_2}^2 = \lim_{t\rightarrow \infty} \mathbb{E}[z^{\sf T}(t)z(t)]\,,
\end{equation}
when each component of $\eta(t)$ is stochastic white noise with unit covariance (i.e., $\mathbb{E}[\eta(t)\,\eta(t^\prime)^{\sf T}]=\delta(t-t^\prime)I$). Therefore, $\|G\|_{\mathcal{H}_2}$ measures how much the output varies in steady-state under stochastic disturbances.

If the state matrix $A$ is Hurwitz, then the $\mathcal{H}_2$ norm is finite, and can be computed as \cite[Ch.4]{KZ-JCD:98}
\begin{equation}\label{Eq:H2Obsv}
\|G\|_{\mathcal{H}_2}^2 = \mathrm{Tr}(B^{\sf T}XB)\,,
\end{equation}
where the observability Gramian $X = X^{\sf T} \succeq \vzeros[]$ is the unique solution to the Lyapunov equation
\begin{equation}\label{Eq:ObsvGram}
 A^{\sf T}X+ XA + C^{\sf T}C = \vzeros[].
\end{equation}
If the pair $(C,A)$ is observable, then $X$ is positive-definite.

\section{$\mathcal{H}_2$ Performance of Saddle-Point Methods}
\label{Section: Performance}

We now subject the saddle-point dynamics \eqref{Eq:SaddlePointContinuous} to disturbances in both the primal and dual equations. Specifically, we assume the vectors $b$ and $c$ are subject to disturbances $\eta_{b} \in \real^{n_b}$ and $\eta_c \in \real^{n_x}$, and make the substitutions $b \mapsto b + t_b\eta_{b}$ and $c \mapsto c + t_c\eta_{\rm c}$ in the saddle point dynamics \eqref{Eq:SaddlePointContinuous}. The scalar parameters $t_b, t_c \geq 0$ parameterize the strength of the disturbances, and will help us keep track of which terms in the resulting norm expressions arise from which disturbances. As an example, when we study distributed resource allocation problems in Section \ref{Sec:Resource}, $b$ will have the interpretation of a vector of demands for some resource, and $\eta_{b}$ will therefore model a fluctuation or disturbance to this nominal demand.

After translating the nominal equilibrium point \eqref{Eq:Equilibrium} of the system to the origin\footnote{In the remainder of the section we assume that we have made the change of state variables $\Delta x = x - x^{\star}$, $\Delta \nu = \nu - \nu^{\star}$, and with an abuse of notation we drop the $\Delta$'s and simply refer to the error coordinates as $x$ and $\nu$.}, we obtain the LTI system

\begin{subequations}\label{Eq:SaddlePointIO}
\begin{align}
	\begin{bmatrix}
		{\mathcal{T}_x}\dot{x}\\
		{\mathcal{T}_\nu}\dot{\nu}
	\end{bmatrix}
&=
	\begin{bmatrix}
		-Q & -S^{\sf T}\\
		S & \vzeros[]
	\end{bmatrix}
	\begin{bmatrix}
		x\\
		\nu
	\end{bmatrix}
	-
	\begin{bmatrix}t_c I_{n_x} & \vzeros[]\\ \vzeros[] & t_bW_b\end{bmatrix}\begin{bmatrix}
	\eta_{c} \\ \eta_{b}
	\end{bmatrix}\,, \\
z &= \begin{bmatrix}C_1 & \vzeros[] \end{bmatrix}\begin{bmatrix}x \\ \nu \end{bmatrix}\,,
\end{align}
\end{subequations}
where $C_1 \in \real^{n_x \times n_x}$ is an output matrix.

As the system \eqref{Eq:SaddlePointIO} is written in error coordinates, convergence to the saddle-point optimizer $(x^{\star},\nu^{\star})$ from  \eqref{Eq:Equilibrium} is equivalent to convergence of $(x(t),\nu(t))$ to the origin. {How should we measure this convergence? A natural way is to use the cost matrix $Q$ from the optimization problem \eqref{Eq:ConvexEqualityConstrained} as a weighting matrix, and to study the performance output $\|z(t)\|_2^2 = x(t)^{\sf T}Qx(t)$, which is obtained by choosing $C_1 = Q^{\frac{1}{2}}$. For example, in the context of the resource allocation problems in Section \ref{Sec:Resource}, these weights describe the relative importance of the various resources.} 

\bigskip

\begin{theorem}[\bf Saddle-Point Performance]\label{Thm:PrimalOutputs}
Consider the input-output saddle-point dynamics \eqref{Eq:SaddlePointIO} with $Q = Q^{\sf T} \succ \vzeros[]$ diagonal, and let $C_1 = Q^{\frac{1}{2}}$ so that $\|z(t)\|_2^2 = x^{\sf T}(t)Qx(t)$. Then the squared $\mathcal{H}_2$ norm of the saddle-point system \eqref{Eq:SaddlePointIO} is
\begin{equation}
\begin{aligned}\label{Eq:H2CostOutput1}
\|G\|_{\mathcal{H}_2}^2 &= \frac{t_c^2}{2}\mathrm{Tr}({\mathcal{T}_x^{-1}}) + \frac{t_b^2}{2}\mathrm{Tr}(W_b^{\sf T}{\mathcal{T}_\nu^{-1}}W_b)\,.
\end{aligned}
\end{equation}
\end{theorem}

\begin{pfof}{Theorem \ref{Thm:PrimalOutputs}}
We will directly construct the unique positive-definite observability Gramian; since the system is internally stable (Lemma \ref{Lem:Stable}), this also indirectly establishes observability \cite[Exercise 4.8.1]{GED-FP:00}. Assuming for the moment a block-diagonal observability {Gramian} $X = \mathrm{blkdiag}(X_1,X_2)$, the Lyapunov equation \eqref{Eq:ObsvGram} yields the two equations

\begin{subequations}
\label{Eq:ObsGramSub}
\begin{align}
\label{Eq:ObsGramSub1}
X_1{\mathcal{T}_x^{-1}}Q + Q{\mathcal{T}_x^{-1}}X_1 - C_1^2 &= \vzeros[]\,,\\
\label{Eq:ObsGramSub2}
X_2{\mathcal{T}_\nu^{-1}}S - S{\mathcal{T}_x^{-1}}X_1 &= \vzeros[]\,,
\end{align}
\end{subequations}
with the third independent equation trivially being $\vzeros[] = \vzeros[]$. By inspection, the solution to \eqref{Eq:ObsGramSub1} is diagonal and given by
$$
X_1 = \frac{1}{2}{\mathcal{T}_x}\,Q^{-1} C_1^2 = \frac{1}{2}{\mathcal{T}_x},
$$
since $Q$ is diagonal and $C_1 = Q^{\frac{1}{2}}$.
Clearly $X_1$ is positive definite and symmetric. Since $S$ has full row-rank, and $X_2$ can be uniquely recovered from \eqref{Eq:ObsGramSub2} as
$$
X_2 = S{\mathcal{T}_x^{-1}}X_1S^{\dagger}{\mathcal{T}_\nu} = \frac{1}{2}SS^{\dagger}{\mathcal{T}_\nu} = \frac{1}{2}{\mathcal{T}_\nu}\,.
$$
It follows that $X_2$ is positive definite, and therefore $X = \frac{1}{2}\mathrm{blkdiag}({\mathcal{T}_x,\mathcal{T}_\nu})$ is the unique positive definite solution to \eqref{Eq:ObsvGram}. Since $X$ is block diagonal, we find from \eqref{Eq:H2Obsv} and \eqref{Eq:SaddlePointIO}
\begin{align*}
\|G\|_{\mathcal{H}_2}^2 &= t_c^2\,\mathrm{Tr}({\mathcal{T}_x^{-1}}X_{1}{\mathcal{T}_x^{-1}}) + t_b^2\,\mathrm{Tr}(W_b^{\sf T}{\mathcal{T}_\nu^{-1}}X_2{\mathcal{T}_\nu^{-1}}W_b)\,,
\end{align*}
from which the result follows.
\end{pfof}

\medskip
We make two key observations about the result \eqref{Eq:H2CostOutput1}. First, \eqref{Eq:H2CostOutput1} is independent of both the cost matrix $Q$ and the constraint matrix $S$; neither matrix has any influence on the value of the system norm. Second, the expression in \eqref{Eq:H2CostOutput1} scales inversely with the time constants {$\mathcal{T}_x$} and {$\mathcal{T}_\nu$}, which indicates an inherent trade-off between convergence speed and input-output performance. As a special case of Theorem \ref{Thm:PrimalOutputs}, suppose that {$\mathcal{T}_x, \mathcal{T}_\nu$} are multiples of the identity matrix, i.e., ${\mathcal{T}_x} = \bar{\tau}_{x} I_{n_x}$, ${\mathcal{T}_\nu} = \bar{\tau}_{\nu} I_{n_r}$, and that $W_b = I_{n_r}$, meaning there is one disturbance for each constraint. Then \eqref{Eq:H2CostOutput1} reduces to
\begin{equation}\label{Eq:H2CostOutput2}
\|G\|_{\mathcal{H}_2}^2 = \frac{t_c^2}{2\bar{\tau}_{x}}n_x + \frac{t_b^2}{2\bar{\tau}_{\nu}}n_r\,,
\end{equation}
meaning the squared $\mathcal{H}_2$ norm scales linearly in the number of disturbances to the primal dynamics and the number of disturbances to the dual dynamics. While this scaling is quite reasonable, the lack of tuneable controller gains other than the time constants means that convergence speed and input-output performance are always conflicting objectives. Finally, we note that \eqref{Eq:H2CostOutput2} can be immediately reinterpreted as a design equation for the time constants. That is, given a specified level $\gamma > 0$ of desired $\mathcal{H}_{2}$ performance, one has $\|G\|_{\mathcal{H}_2} \leq \gamma$ if
$$
\min\{\bar{\tau}_x,\bar{\tau}_{\nu}\} \geq \frac{1}{\gamma^2}\left(\frac{t_c^2 n_x}{2} + \frac{t_b^2 n_r}{2}\right).
$$

\subsection{Performance of Regularized Saddle-Point Methods}
\label{Sec:Regularized}

A common variation of the Lagrangian optimization framework includes a quadratic penalty \cite{MBK-NL:16,AS-GL:14,ED-AS:16} on the dual variable $\nu$ in the Lagrangian \eqref{Eq:Lagrangian}. The so-called \emph{regularized} Lagrangian assumes the form
\begin{equation}\label{Eq:Lagrangian_reg}
L_{\rm reg}(x,\nu) = \frac{1}{2}x^{\sf T}Qx + c^{\sf T}x + \nu^{\sf T}(Sx-W_bb)-\frac{\epsilon}{2} \|\nu\|_2^2\,,
\end{equation}
where $\epsilon > 0$ is small. The regularization term adds concavity to the Lagrangian, and has been shown to increase the convergence rate of optimization algorithms. However, this regularization alters the equilibrium of the closed-loop system, which moves from the value in \eqref{Eq:Equilibrium} to the new value
\begin{equation}\label{Eq:Equilibrium_reg}
\begin{bmatrix}
x_{\rm reg}^{\star} \\ \nu_{\rm reg}^{\star}
\end{bmatrix} = 
\begin{bmatrix}
-Q^{-1}(S^{\sf T}\nu_{\rm reg}^{\star} + c )\\
-(SQ^{-1}S^{\sf T}+\epsilon I_{n_r})^{-1}(W_bb+SQ^{-1}c)
\end{bmatrix}.
\end{equation}
The penalty coefficient $\epsilon$ is chosen to strike a balance between the convergence rate improvement and the deviation of $(x_{\rm reg}^{\star},\nu_{\rm reg}^{\star})$ from $(x^{\star},\nu^{\star})$. The Lagrangian \eqref{Eq:Equilibrium_reg} also admits a continuous-time saddle-point algorithm with {\it regularized saddle-point dynamics}
\begin{equation}
\label{Eq:SaddlePointContinuous_reg}
\begin{aligned}
{\mathcal{T}_x}\dot{x} &= -Qx - S^{\sf T}\nu - c\\
{\mathcal{T}_\nu} \dot{\nu} &= Sx-W_b b-\epsilon \nu\,.
\end{aligned}
\end{equation}
Quite strikingly, we shall observe that a small regularization which has a minor effect on the equilibrium, achieves a tremendous improvement in performance.

As we did with the standard saddle-point dynamics, we can shift the equilibrium point $(x^{\star}_{\rm reg},\nu^{\star}_{\rm reg})$ of \eqref{Eq:SaddlePointContinuous_reg} to the origin and introduce disturbances to the parameters $b$ and $c$, leading to the input-output model
\begin{equation}\label{Eq:SaddlePointIO_reg}
\begin{aligned}
	\begin{bmatrix}
		{\mathcal{T}_x}\dot{x}\\
		{\mathcal{T}_\nu}\dot{\nu}
	\end{bmatrix}
&=
	\begin{bmatrix}
		-Q & -S^{\sf T}\\
		S & -\epsilon I_{n_r}
	\end{bmatrix}
	\begin{bmatrix}
		x\\
		\nu
	\end{bmatrix}
	-
	\begin{bmatrix}t_cI_{n_x} & \vzeros[]\\ \vzeros[] & t_bW_b\end{bmatrix}\begin{bmatrix}
	\eta_{c} \\ \eta_{b}
	\end{bmatrix}\, \\
z &= \begin{bmatrix}Q^{\frac{1}{2}} & \vzeros[] \end{bmatrix}\begin{bmatrix}x \\ \nu \end{bmatrix}\,,
\end{aligned}
\end{equation}
where we consider the time constant matrices {$\mathcal{T}_x$, $\mathcal{T}_\nu$}, and disturbances $\eta_b$, $\eta_c$ as in \eqref{Eq:SaddlePointIO}. \\

\begin{theorem}[\bf Regularized Saddle-Point Performance]\label{Thm:PrimalOutputs_reg}
Consider the input-output regularized saddle-point dynamics \eqref{Eq:SaddlePointIO_reg} denoted by $G_\text{reg}$ with $Q = Q^{\sf T} \succ \vzeros[]$ a diagonal matrix. Then the squared $\mathcal{H}_2$ norm $\|G_\text{reg}\|_{\mathcal{H}_2}^2$ of the system \eqref{Eq:SaddlePointIO_reg} is upper-bounded by \eqref{Eq:H2CostOutput1} with strict inequality. 
\end{theorem}

\begin{pfof}{Theorem \ref{Thm:PrimalOutputs_reg}}
We rewrite \eqref{Eq:SaddlePointIO_reg} in the standard state-space form $G_{\rm reg} := (A_{\rm reg}, B, C,\vzeros[])$, where
\begin{subequations}\label{Eq:SaddlePointIO_reg_stan}
\begin{align*}
	\begin{bmatrix}
		\dot{x}\\
		\dot{\nu}
	\end{bmatrix}
&=
	\underbrace{\begin{bmatrix}
		-{\mathcal{T}_x}^{-1}Q & -{\mathcal{T}_x}^{-1}S^{\sf T}\\
		{\mathcal{T}_\nu}^{-1}S & -{\mathcal{T}_\nu}^{-1}\epsilon I_{n_r}
	\end{bmatrix}}_{A_\text{reg}}
	\begin{bmatrix}
		x\\
		\nu
	\end{bmatrix}\\
	&\quad -
	\underbrace{\begin{bmatrix}{\mathcal{T}_x}^{-1} t_c & \vzeros[]\\ \vzeros[] & {\mathcal{T}_\nu}^{-1}W_b t_b\end{bmatrix}}_{B}\begin{bmatrix}
	\eta_{c} \\ \eta_{b}
	\end{bmatrix}\,,
z = \underbrace{\begin{bmatrix}Q^{\frac{1}{2}} & \vzeros[] \end{bmatrix}}_{C}\begin{bmatrix}x \\ \nu \end{bmatrix}\,.
\end{align*}
\end{subequations}
\end{pfof}
One may verify that $A_{\rm reg}$ is Hurwitz and that $G_{\rm reg}$ is observable. Consider the observability Gramian from \eqref{Eq:ObsGramSub}, i.e., $X = \frac{1}{2}\mathrm{blkdiag}({\mathcal{T}_x,\mathcal{T}_{\nu}})$. An easy computation shows that
\begin{equation}
\label{e:regGramian}
X\, A_\text{reg}+A_\text{reg}^{\sf T}X+C^{\sf T} C= \begin{bmatrix}\vzeros[] & \vzeros[]\\ \vzeros[] & -\epsilon I_{n_r}\end{bmatrix} \preceq \vzeros [].
\end{equation}
We conclude that $X$ is a \emph{generalized} observability Gramian for the regularized system $G_\text{reg}$. If $X_{\epsilon}$ is the true observability Gramian for $G_{\rm reg}$, then $X_{\epsilon} \neq X$ and $X_{\epsilon} \preceq X$ \cite[Chapter 4.7]{GED-FP:00}, and we conclude that
\begin{equation}\label{Eq:regInequality}
\|G_\text{reg}\|_{\mathcal{H}_2}^2 = \mathrm{Tr}(B^{\sf T}X_{\epsilon}B) \leq \mathrm{Tr}(B^{\sf T}XB) = \|G\|_{\mathcal{H}_{2}}^2\,.
\end{equation}
It remains only to show that the above inequality holds strictly. Proceeding by contradiction, assume that $\mathrm{Tr}(B^{\sf T}X_{\epsilon}B) = \mathrm{Tr}(B^{\sf T}XB)$, which implies that $\mathrm{Tr}(B^{\sf T}(X-X_{\epsilon})B) = 0$. Since $X - X_{\epsilon} \succeq \vzeros[]$, we may write $X - X_{\epsilon} = F^{\sf T}F$ for some matrix $F$, and
$$
0 = \mathrm{Tr}(B^{\sf T}(X-X_{\epsilon})B) = \mathrm{Tr}(B^{\sf T}F^{\sf T}FB)\,.
$$
Since $B$ has full row rank, this implies that $F$ must be zero, and thus $X = X_{\epsilon}$ which is a contradiction.
\begin{corollary}[\bf Regularized Saddle-Point Performance with One Constraint]\label{Prop:performance-gap}
Consider the case with one constraint ($n_r = 1$) and one constraint disturbance ($n_b = 1$) with uniform problem parameters $Q = q I_{n_x}$, ${\mathcal{T}_x} = \bar{\tau}_x I_{n_x}$, ${\mathcal{T}_\nu} = \bar \tau_\nu I_{n_r} $ for scalars $q, \bar \tau_x, \bar \tau_\nu > 0$ and $W_b = 1$. Then, we have
\begin{equation}\label{Eq:GregExpression}
\|G\|_{\mathcal{H}_2}^2 -\|G_\text{reg}\|_{\mathcal{H}_2}^2= \alpha_\epsilon  t_c^2 + \gamma_\epsilon {t_b^2} , 
\end{equation}
where $s = \|S\|_2$ and
\[
\alpha_\epsilon=\frac{\epsilon s^2}{2(\epsilon q + s^2)(\epsilon\bar \tau_x  + q\bar \tau_\nu)} \;,
\]
\[
\gamma_\epsilon = \frac{\epsilon (\bar \tau_x q  \epsilon + q^2 \bar \tau_\nu  + \bar \tau_xs^2 )}{2\bar \tau_\nu(\epsilon q + s^2)(\epsilon \bar \tau_x  + q \bar \tau_v)}\,.
\]
\end{corollary} 

\medskip{}
\begin{pfof}{Corollary \ref{Prop:performance-gap}}
Let $\Delta:=X-X_\epsilon$, where $X = \frac{1}{2}\mathrm{blkdiag}({\mathcal{T}_x,\mathcal{T}_{\nu}})$ is the observability Gramian of \eqref{Eq:SaddlePointIO} with $C_1=Q^{\frac{1}{2}}$, and $X_\epsilon$ is the observability Gramian of $G_\text{reg}$. As noted in the proof of Theorem \ref{Thm:PrimalOutputs_reg}, $X\succeq X_\epsilon$, and thus the matrix $\Delta$ is positive semidefinite. Clearly, $\Delta$ satisfies
\[
A_\text{reg}^{\sf T} (X-\Delta) + (X-\Delta) A_\text{reg} +C^{\sf T}C=0.
\]
By \eqref{e:regGramian}, this reduces to 
\begin{equation}\label{Eq:eps-sys-obs}
A_\text{reg}^{\sf T} \Delta + A_\text{reg} \Delta +  \begin{bmatrix}\vzeros[] & \vzeros[]\\ \vzeros[] & \epsilon I_{n_r}\end{bmatrix} = \vzeros[]\,.
\end{equation}
Then it is easy to see that
\begin{equation}\label{Eq:improved-H2}
\|G_\text{reg}\|_{\mathcal{H}_2}^2= \|G\|_{\mathcal{H}_2}^2- \|G_\epsilon\|_{\mathcal{H}_2}^2,
\end{equation}
where $\|G\|_{\mathcal{H}_2}^2$ is as in \eqref{Eq:H2CostOutput1}, and $G_\epsilon(s):=C_\epsilon (sI- A_\text{reg})^{-1} B$ with $C_\epsilon= \begin{bmatrix} \vzeros[] & \epsilon I_{n_r}\end{bmatrix}$. Therefore, the improvement in the $\mathcal{H}_2$-norm performance is  equal to the squared $\mathcal{H}_2$-norm of the axillary system given by $G_\epsilon$. 
Next, we calculate the $\mathcal{H}_2$-norm of $G_\epsilon$, which requires computing the observability Gramian from 
\eqref{Eq:eps-sys-obs}. Consider the matrix
\begin{equation}\label{Eq:X-eps}
{\bar \Delta}=\begin{bmatrix}
\alpha S^{\sf T}S & \beta S^{\sf T}\\
\beta S & \gamma 
\end{bmatrix}
\end{equation}
where $\alpha$, $\beta$, $\gamma$ are constant and positive.
Straightforward calculation shows that by choosing 
\begin{equation}\label{Eq:abg}
\begin{bmatrix}
\alpha\\
\beta\\
\gamma\\
\end{bmatrix}
=
\begin{bmatrix}
\frac{s^2}{\bar \tau_x} & \frac{\epsilon}{\bar \tau_\nu}+\frac{q}{\bar \tau_x} & -\frac{1}{\bar \tau_\nu}\\[1.5mm]
\frac{q}{\bar \tau_x}  & -\frac{1}{\bar \tau_\nu} & 0\\[1.5mm]
0 & \frac{\bar \tau_\nu s^2}{\bar \tau_x \epsilon} & 1
\end{bmatrix}
^{-1}
\begin{bmatrix}
 0\\[1.5mm]
 0\\[1.5mm]
\frac{ \bar \tau_v}{2}
\end{bmatrix},
\end{equation}
the matrix $\bar \Delta$ is a solution to the Lyapunov equation {\eqref{Eq:eps-sys-obs}}. Given the fact that $A_\text{reg}$ is Hurwitz,  
this solution is unique and the matrix $\Delta=\overline \Delta$ is the observability Gramian of the system given by $G_\epsilon$.
The proof is completed by calculating the inverse in \eqref{Eq:abg}, using \eqref{Eq:H2Obsv}, and noting $\mathrm{Tr} (S^{\sf T}S)=s^2$.
\end{pfof}

\begin{figure*}[!ht]
        \begin{center}
        \begin{subfigure}[!ht]{\columnwidth}
        		\centering
                \includegraphics[height=0.75\columnwidth]{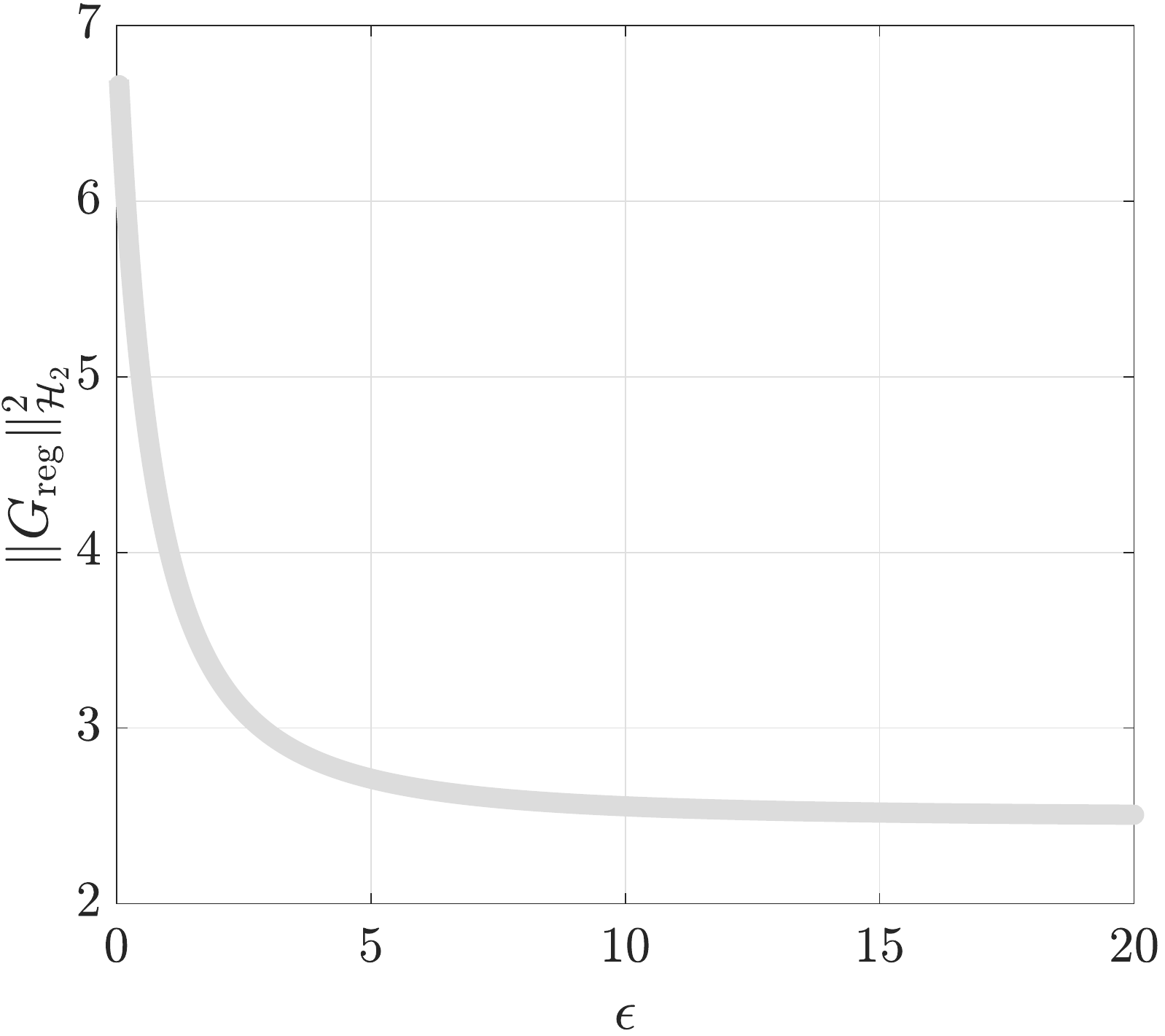}
                \caption[]{q=3}
                \label{Fig:RA4}
        \end{subfigure}%
        \begin{subfigure}[!ht]{\columnwidth}
        		\centering
                \includegraphics[height=0.75\columnwidth]{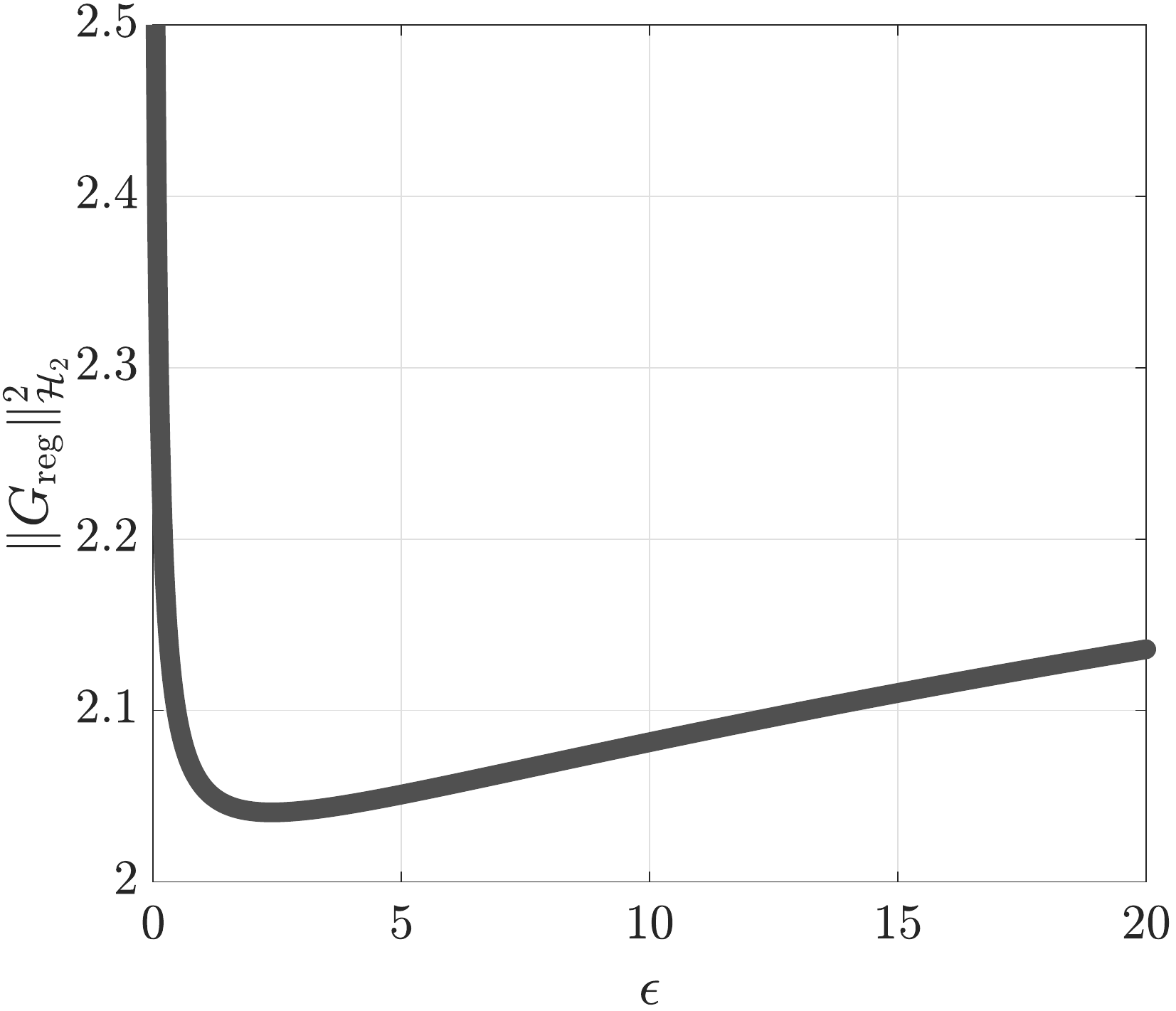}
                       \caption[]{q=0.05}
                \label{Fig:RA}
                \end{subfigure}
		\caption[]{System norm of regularized dynamics as a function of $\epsilon$, for parameters $\bar{\tau}_x = \bar{\tau}_{\nu} = t_c = t_b = 1$, $Q = q I_{5}$, and $S =[0.82\,\,    0.90\,\,    0.13\,\,    0.91\,\,    0.63] $.}
		\label{Fig:Regularization}
\end{center}
\end{figure*}

Corollary \ref{Prop:performance-gap} quantifies the performance improvement resulting  from the regularization in the special case of a single constraint and uniform parameters.  For sufficiently large $\epsilon$,  this improvement is approximated by $\frac{t_b^2}{2\bar{\tau}_{\nu}}$, which
coincides with the second term in the right-hand side of \eqref{Eq:H2CostOutput2}. This means that, as expected, the constraints do not contribute to the $\mathcal{H}_2$-norm in case the penalty term in the regularized Lagrangian \eqref{Eq:Lagrangian_reg} tends to infinity.
On the other hand, for $\epsilon \rightarrow 0^+$, the $\mathcal{H}_2$ norm of the regularized dynamics \eqref{Eq:SaddlePointIO_reg} clearly converges to the $\mathcal{H}_2$ norm of \eqref{Eq:SaddlePointIO}. In general, the improvement in the input-output performance obtained due to regularization is not a monotonic function of the regularization parameter $\epsilon$. To illustrate this, Figure \ref{Fig:Regularization} plots the right-hand side of \eqref{Eq:GregExpression} as a function of $\epsilon$. It is noteworthy that for both the plots, even a modest $\epsilon$ improves the performance. Depending on the specific values of the parameters, the maximum performance gain may be achieved as $\epsilon \rightarrow \infty$ (Figure \ref{Fig:Regularization}(a)) or at a finite value of $\epsilon$ (Figure \ref{Fig:Regularization}(b)).

\subsection{Performance of Augmented Saddle-Point Methods}
\label{Sec:Augmented}
Another option for improving the $\mathcal{H}_2$ performance of saddle-point methods is to return to the Lagrangian function \eqref{Eq:Lagrangian} and instead consider the \emph{augmented Lagrangian} \cite{DPB:08}
\begin{equation}\label{Eq:LAug}
L_{\text{aug}}(x,\nu) \triangleq L(x,\nu) + \frac{\rho}{2}\|Sx-W_b b\|_2^{2}\,,
\end{equation}
where we have incorporated the squared constraint $Sx - W_b b = \vzeros[n_r]$ into the Lagrangian with a gain $\rho \geq 0$. One way to interpret this is that the term $\frac{\rho}{2}\|Sx-W_b b\|_2^{2}$ adds additional convexity to the Lagrangian in the $x$ variable.

It follows that $(x,\nu)$ is a saddle point of $L_{\text{aug}}(x,\nu)$ if and only if it is a saddle point of $L(x,\nu)$, and hence the optimizer is unchanged. Applying the saddle-point method to the augmented Lagrangian $L_{\text{aug}}(x,\nu)$, we obtain the \emph{augmented saddle-point dynamics}
\begin{equation}\label{Eq:SaddlePointAugmented}
\begin{aligned}
{\mathcal{T}_x}\dot{x} &= -(Q+\rho S^{\sf T}S)x - S^{\sf T}\nu - c + \rho S^{\sf T}W_b b\\
{\mathcal{T}_\nu} \dot{\nu} &= Sx-W_b b\,.
\end{aligned}
\end{equation}
One may verify that as before, the unique stable equilibrium point of \eqref{Eq:SaddlePointAugmented} is given by \eqref{Eq:Equilibrium}. We again consider disturbances $\eta_{b}$ and $\eta_c$, and make the substitution $b \mapsto b + t_b\eta_{b}$ and $c \mapsto c + t_c\eta_{\rm c}$. After translating the equilibrium point to the origin, we obtain the MIMO system
\begin{equation}\label{Eq:SaddlePointAugmentedIO}
\begin{aligned}
	\begin{bmatrix}
		{\mathcal{T}_x}\dot{x}\\
		{\mathcal{T}_\nu}\dot{\nu}
	\end{bmatrix}
&=
	\begin{bmatrix}
		-(Q+\rho S^{\sf T}S) & -S^{\sf T}\\
		S & \vzeros[]
	\end{bmatrix}
	\begin{bmatrix}
		x\\
		\nu
	\end{bmatrix}\\
	&\qquad \qquad \qquad
	-
	\begin{bmatrix}
	t_cI_{n_x} & -\rho t_b S^{\sf T}W_b\\  \vzeros[] & t_bW_b
	\end{bmatrix}\begin{bmatrix}
	\eta_{c} \\ \eta_{b}
	\end{bmatrix}\\
z &= \begin{bmatrix}Q^{\frac{1}{2}} & \vzeros[] \end{bmatrix}\begin{bmatrix}x \\ \nu \end{bmatrix}\,.
\end{aligned}
\end{equation}
The additional term $-\rho S^{\sf T}S$ in the dynamics \eqref{Eq:SaddlePointAugmentedIO} complicates the solution of the Lyapunov equation, and we require additional assumptions to obtain an explicit formula. We consider the parametrically uniform case where $Q = q I_{n_x}$, ${\mathcal{T}_x} = \bar{\tau}_{x} I_{n_x}$, and ${\mathcal{T}_\nu} = \bar{\tau}_{\nu} I_{n_r}$  for scalars $q, \bar{\tau}_{x}, \bar{\tau}_{\nu} > 0$. The next result appeared in \cite{JWSP-BKP-NM-FD:16c} without proof.

\medskip
\begin{theorem}[\bf {Augmented Saddle-Point Performance with Uniform Parameters}]\label{Thm:AugmentedPrimalOutputs}
Consider the input-output augmented saddle-point dynamics \eqref{Eq:SaddlePointAugmentedIO}, denoted by $G_\text{aug}$ under the above assumptions, with performance output $\|z(t)\|_2^2 = x^{\sf T}(t)Qx(t) = q\|x(t)\|_2^2$. Then the squared $\mathcal{H}_2$ norm of the augmented saddle-point system \eqref{Eq:SaddlePointAugmentedIO} for identically weighted disturbances $W_b=I_{n_r}$ is
\begin{equation}
\begin{aligned}\label{Eq:H2CostOutputAugmented1}
\|G_\text{aug}\|_{\mathcal{H}_2}^2 &= \frac{t_c^2}{2\bar{\tau}_{x}}(n_x-n_r) + \left(\frac{t_b^2}{2\bar{\tau}_{\nu}}+\frac{t_c^2}{2\bar{\tau}_{x}}\right)\sum_{i=1}^{n_r}\frac{q}{q+\rho \sigma_i^2}\\
&\quad + \frac{t_b^2}{2\bar{\tau}_{x}}\sum_{i=1}^{n_r}\frac{q\rho^2\sigma_i^2}{q+\rho \sigma_i^2}\,,
\end{aligned}
\end{equation}
where $\sigma_i$ is the $i$th non-zero singular value of $S$.
\end{theorem}

\begin{pfof}{Theorem \ref{Thm:AugmentedPrimalOutputs}}
Under the given assumptions, the system \eqref{Eq:SaddlePointAugmentedIO} simplifies to
\begin{subequations}
\renewcommand*{\arraystretch}{1.2}
\begin{align*}
	\begin{bmatrix}
		\dot{x}\\
		\dot{\nu}
	\end{bmatrix}
&=
	\begin{bmatrix}
		-\frac{1}{\bar{\tau}_{x}}(qI_{n_x}+\rho S^{\sf T}S) & -\frac{1}{\bar{\tau}_{x}}S^{\sf T}\\
		\frac{1}{\bar{\tau}_{\nu}}S & \vzeros[]
	\end{bmatrix}
	\begin{bmatrix}
		x\\
		\nu
	\end{bmatrix}\\
	&\quad -
	\begin{bmatrix}\frac{t_c}{\bar{\tau}_{x}}I_{n_x} & -\frac{t_b\rho}{\bar{\tau}_{x}}S^{\sf T}{W_b}\\ \vzeros[] &{\frac{t_b}{\bar{\tau}_{\nu}}{W_b}} \end{bmatrix}\begin{bmatrix}
	\eta_{c} \\ \eta_{b}
	\end{bmatrix}\,,\,\,
z = \begin{bmatrix}q^{\frac{1}{2}}I_{n_x} & \vzeros[] \end{bmatrix}\begin{bmatrix}x \\ \nu \end{bmatrix}\,.
\end{align*}
\end{subequations}
Let $S = U\Sigma V^{\sf T}$ be the singular value decomposition of $S$, where $U \in \real^{n_r \times n_r}$ and $V \in \real^{n_x \times n_x}$ are both orthogonal matrices, $\Sigma = \begin{bmatrix}\bar{\Sigma} &\vzeros[n_r\times(n_x-n_r)]\end{bmatrix}$ and $\bar{\Sigma} \in \real^{n_r \times n_r}$ is the diagonal matrix of non-zero singular values. Consider now the invertible change of variables $\tilde{x} = V^{\sf T}x$, $\tilde{\nu} = U^{\sf T}\nu$. In these new coordinates, the dynamics become
\begin{subequations}
\renewcommand*{\arraystretch}{1.2}
\begin{align*}
	{\begin{bmatrix}
		\dot{\tilde{x}}\\
		\dot{\tilde{\nu}}
	\end{bmatrix}}
&=
	\begin{bmatrix}
		-\frac{1}{\bar{\tau}_{x}}(qI_{n_x}+\rho \Sigma^{\sf T}\Sigma) & -\frac{1}{\bar{\tau}_{x}}\Sigma^{\sf T}\\
		\frac{1}{\bar{\tau}_{\nu}}\Sigma & \vzeros[]
	\end{bmatrix}
	\begin{bmatrix}
		\tilde{x}\\
		\tilde{\nu}
	\end{bmatrix}\\
	&-
	\begin{bmatrix}{\frac{t_c}{\bar{\tau}_{x}}V^{\sf T}} & -\frac{t_b\rho}{\bar{\tau}_{x}}V^{\sf T}S^{\sf T}{W_b} \\ \vzeros[] & {\frac{t_b}{\bar{\tau}_{\nu}}}{U^{\sf T} {W_b}}\end{bmatrix}\begin{bmatrix}
	\eta_{c} \\ \eta_{b}
	\end{bmatrix}\,,\,\,
{{y}} = \begin{bmatrix}q^{\frac{1}{2}}{V} & \vzeros[] \end{bmatrix}\begin{bmatrix}\tilde{x} \\ \tilde{\nu} \end{bmatrix}\,.
\end{align*}
\end{subequations}
We now show the observability of the pair $(C,A)$; note that this is equivalent to observability of $(C^{\sf T}C, A)$. First note that $\mathrm{ker}(C^{\sf T}C)$ is spanned by $\begin{bmatrix}\vzeros[n_x]^{\sf T} & \nu^{\sf T}\end{bmatrix}^{\sf T}$. Now, suppose that $\begin{bmatrix}\vzeros[n_x]^{\sf T} & \nu^{\sf T}\end{bmatrix}^{\sf T}$ is an eigenvector of $A$ with eigenvalue $\lambda$:
\begin{equation*}
\begin{bmatrix}
		-\frac{1}{\bar{\tau}_{x}}(qI_{n_x}+\rho \Sigma^{\sf T}\Sigma) & -\frac{1}{\bar{\tau}_{x}}\Sigma^{\sf T}\\
		\frac{1}{\bar{\tau}_{\nu}}\Sigma & \vzeros[]
	\end{bmatrix}
	\begin{bmatrix}
		\vzeros[n_x]\\
		{\nu}
	\end{bmatrix}= \lambda \begin{bmatrix}
		\vzeros[n_x]\\
		{\nu}
	\end{bmatrix}.
\end{equation*}
Since by stability $\mathrm{Re}(\lambda) <0$, the above relation only holds for $\nu=\vzeros[n_r]$, which shows observability by the eigenvector test. Assuming a block-diagonal observability {Gramian} $X = \mathrm{blkdiag}(X_1,X_2)$, the Lyapunov equation \eqref{Eq:ObsvGram} yields the two independent equations
\begin{subequations}
\begin{align}
\label{Eq:ObsGramNew1} 
X_1\frac{1}{\bar{\tau}_{x}}(qI_{n_x}+\rho \Sigma^{\sf T}\Sigma) + \frac{1}{\bar{\tau}_{x}}(qI_{n_x}+\rho \Sigma^{\sf T}\Sigma)X_1 &= q I_{n_x} \,,\\
\label{Eq:ObsGramNew2}
X_2\frac{1}{\bar{\tau}_{\nu}}\Sigma - \Sigma\frac{1}{\bar{\tau}_{x}}X_1 &= \vzeros[]\,,
\end{align}
\end{subequations}
where we have used the fact that $V^{\sf T}V = I_{n_x}$. Noting that $\Sigma^{\sf T}\Sigma = \mathrm{blkdiag}(\bar{\Sigma},\vzeros[(n_x-n_r)\times(n_x-n_r)])$, we find by inspection that the solution to \eqref{Eq:ObsGramNew1} is diagonal and given by
\begin{align*}
{X_1}_{ii} &= \dfrac{1}{2} \dfrac{q}{q+\rho \sigma_i^2}\bar{\tau}_{x}\,,& \quad i \in \{&1,\ldots,n_r\}\,,\\
{X_1}_{ii} &= \dfrac{1}{2}\bar{\tau}_{x}\,,& \quad i \in \{&n_r+1, \ldots ,n_x\}\,.
\end{align*}
Observe that $X_1$ is positive definite and symmetric. The matrix equation \eqref{Eq:ObsGramNew2} admits a solution $X_2$ if and only if $\mathrm{ker}(\Sigma) \subseteq \mathrm{ker}(\Sigma X_1)$, which holds in this case since $X_1$ is diagonal. The lower block $X_2$ can therefore be uniquely recovered from \eqref{Eq:ObsGramNew2} as
\begin{align*}
X_2 &= \bar{\tau}_{x}^{-1}(\Sigma X_1\Sigma^{\sf T})(\Sigma \Sigma^{\sf T})^{-1}\bar{\tau}_{\nu} \,.
\end{align*}
A straightforward calculation shows that this is equivalent to the component formula
\begin{align*}
{X_2}_{ii} &= \dfrac{1}{2} \dfrac{q}{q+\rho \sigma_i^2}\bar{\tau}_{\nu}\,, \quad i \in \{1,\ldots,n_r\}\,.
\end{align*}
It follows that $X_2$ is diagonal and positive definite, and therefore $X = \mathrm{blkdiag}(X_1, X_2)$ is the unique positive definite solution. A calculation using \eqref{Eq:H2Obsv} now shows that
\begin{equation*}
\begin{aligned}
\|G_\text{aug}\|_{\mathcal{H}_2}^2 &= \frac{t_c^2}{\bar{\tau}_{x}^2}\mathrm{Tr}\left(V X_{1}V^{\sf T}\right) + \frac{t_b^2}{\bar{\tau}_{\nu}^2}\mathrm{Tr}\left({W_b^{\sf T}}U X_{2}U^{\sf T}{W_b}\right)\\
&\quad + \frac{t_b^2\rho^2}{\bar{\tau}_{x}^2}\mathrm{Tr}\left({W_b^{\sf T}}SVX_1V^{\sf T}S^{\sf T}{W_b}\right).
\end{aligned}
\end{equation*}
For the special case of one disturbance per constraint, i.e., $W_b = I_{n_r}$, the result \eqref{Eq:H2CostOutputAugmented1} follows by applying the cyclic property of the trace operation.
\end{pfof}
\medskip

Under the assumed restrictions on parameters, Theorem \ref{Thm:AugmentedPrimalOutputs} generalizes Theorem \ref{Thm:PrimalOutputs}, since when $\rho = 0$ the expression \eqref{Eq:H2CostOutputAugmented1} reduces to \eqref{Eq:H2CostOutput2}.
Consider now the dependence of the expression \eqref{Eq:H2CostOutputAugmented1} on the augmentation gain $\rho$. First, in the case when $t_b = 0$ (meaning the vector $b$ is not subject to disturbances), then as $\rho \rightarrow \infty$ the expression \eqref{Eq:H2CostOutputAugmented1} reduces to only the first term: in this case, augmentation unambiguously improves input-output performance. In particular, note that a more favourable scaling than in \eqref{Eq:H2CostOutput1} is achieved when $n_r$ is comparable to $n_x$; see the resource allocation problem in Section~\ref{Sec:Resource}.\footnote{As an observation, we note that even when $S^{\sf T}$ is a sparse matrix, $S^{\sf T}S$ typically will not be, and hence the augmented dynamics \eqref{Eq:SaddlePointAugmented} may not be immediately implementable as a distributed algorithm. A notable exception is when $S$ is the transposed incidence matrix of a sparse graph, which gives $S^{\sf T}S$ as the corresponding sparse graph Laplacian; this will occur in Section~\ref{Sec:Resource}.} On the other hand, if $t_b \neq 0$, then as $\rho$ becomes large, the second term in the expression vanishes, while the third term grows without bound. Therefore, a large augmentation gain will lead to poor input-output performance. This behaviour is explained by examining \eqref{Eq:SaddlePointAugmented}: the vector $b$ enters the primal dynamics multiplied by $\rho$, and hence any noise in $b$ is amplified as $\rho$ grows. To remedy this deficiency in the augmented approach, the next section exploits a dual formulation of the optimization problem \eqref{Eq:ConvexEqualityConstrained}.

\section{Dual and Distributed Dual Methods for Improved Saddle-Point Algorithm Performance}
\label{Sec:DualMethods}

This section develops an approach to overcome the performance issues of augmented Lagrangian methods observed in Theorem \ref{Thm:AugmentedPrimalOutputs} when disturbances enter the constraints. To focus in on these problematic disturbances, in this section we ignore possible disturbances to the vector $c$ and set $t_c = 0$. Section \ref{Sec:CentralizedDual} contains a quick examination of dual ascent, before proceeding to a distributed dual formulation in Section \ref{Sec:DistributedDual}.

\subsection{Centralized Dual Ascent}
\label{Sec:CentralizedDual}
To begin, we return to the Lagrangian \eqref{Eq:Lagrangian} of the optimization problem \eqref{Eq:ConvexEqualityConstrained}, and compute 
$$
x^{\star}(\nu) = \argmin_{x \in \real^n} L(x,\nu) = -Q^{-1}(c + S^{\sf T}\nu)\,.
$$
It follows quickly that the dual function $\Phi(\nu)$ is given by
\begin{equation}\label{Eq:DualFunction}
\begin{aligned}
\Phi(\nu)& = \min_{x \in \real^n} L(x,\nu)\\
&=-\frac{1}{2}\nu^{\sf T}SQ^{-1}S^{\sf T}\nu - \nu^{\sf T}(SQ^{-1}c+W_bb)
- \frac{1}{2}c^{\sf T}Q^{-1}c\,.
\end{aligned}
\end{equation}
With the primal variable eliminated, a possible approach is to simply maximize the dual function $\Phi(\nu)$ via gradient ascent. Introducing disturbance inputs $b \mapsto b + t_b\eta_b$ and performance outputs similar to before, and shifting the unique equilibrium point to the origin, one quickly obtains the input-output dual ascent dynamics
\begin{equation}\label{Eq:DualDynamics}
\begin{aligned}
{\mathcal{T}_{\nu}}\dot{\nu} &= -SQ^{-1}S^{\sf T}\nu - t_bW_b\eta_b\\
z &= -Q^{-\frac{1}{2}}S^{\sf T}\nu\,,
\end{aligned}
\end{equation}
where ${\mathcal{T}_{\nu}} \succ \vzeros[]$. The performance of \eqref{Eq:DualDynamics} is characterized by the following result.

\smallskip

\begin{proposition}[\bf Dual Ascent Performance]\label{Prop:Dual-Ascent}
The $\mathcal{H}_2$ norm of the input-output dual ascent dynamics \eqref{Eq:DualDynamics} is given by
$$
\|G\|^2_{\mathcal{H}_2} = \frac{t_b^2}{2}\mathrm{Tr}(W_b^{\sf T}{\mathcal{T}_{\nu}^{-1}}W_b)
$$
\end{proposition}
\begin{pfof}{Proposition \ref{Prop:Dual-Ascent}}
The Lyapunov equation \eqref{Eq:ObsvGram} for this problem takes the form
$$
-X{\mathcal{T}_{\nu}^{-1}}SQ^{-1}S^{\sf T} - SQ^{-1}S^{\sf T}{\mathcal{T}_{\nu}^{-1}}X + SQ^{-1}S^{\sf T} = \vzeros[]
$$
from which we find the unique solution $X = \frac{1}{2}{\mathcal{T}_{\nu}} \succ \vzeros[]$. With $B = \left[\begin{matrix}t_b {\mathcal{T}_{\nu}^{-1}}W_b\end{matrix}\right]$, the result follows by applying \eqref{Eq:H2Obsv}.
\end{pfof}

Comparing the result of Proposition \ref{Prop:Dual-Ascent} to the unaugmented saddle-point result of Theorem \ref{Thm:PrimalOutputs}, we observe that the terms proportional to $t_b^2$ are identical. Therefore, when considering algorithm performance with disturbances entering the constraints, the primal-dual and pure dual-ascent algorithms achieve identical performance.

\subsection{A Distributed Dual Augmented Lagrangian Method}
\label{Sec:DistributedDual}

Building off the dual function \eqref{Eq:DualFunction}, we now derive a modified augmented Lagrangian algorithm, which can overcome the performance issues posed by disturbances affecting the vector $b$. The particulars of the derivation below are tailored towards distributed solutions, which will be discussed further in Section \ref{Sec:Resource} in the context of distributed resource allocation.
With this application in mind, we will focus in on the case where $n_b = n_x$, so that the $i$th component of the disturbance $b$ can be uniquely associated to the $i$th primal variable $x_i$; this assumption can be relaxed in the derivation below as long as one uniquely assigns components of $b$ and the associated columns of $W_{b}$ to a particular agent. We partition each of the following matrices according to their columns as
$$
\begin{aligned}
S &= \begin{bmatrix}
s_1\,\, s_2 \,\, \cdots \,\, s_{n_x}
\end{bmatrix}, &W_b &= \begin{bmatrix}
w_1 \,\, w_2 \,\,\cdots \,\, w_{n_x}
\end{bmatrix}.
\end{aligned}
$$
With this partitioning, one may quickly see that for $Q=\diag(q_1,\ldots,q_{n_x})$, the dual function \eqref{Eq:DualFunction} may be written as
$$
\Phi(\nu) = \sum_{i=1}^{n_x} \underbrace{\left[-\frac{1}{2q_i}\nu^{\sf T}s_is_i^{\sf T}\nu - \nu^{\sf T}\left(\frac{c_i}{q_i}s_i+ w_ib_i\right) -\frac{c_i^2}{2q_i} \right]}_{:=\widetilde{\Phi}_i(\nu)}.
$$
The dual function appears to separate into a sum, except for the common multiplier $\nu$. To complete the separation, for each $i \in \{1,\ldots,n_x\}$ we introduce a local copy $\nu^i \in \real^{n_r}$ of the vector of Lagrange multipliers $\nu \in \real^{n_r}$, and require that $\nu^i = \nu^j$ for all $i,j \in \{1,\ldots,n_x\}$. To enforce these so-called agreement constraints, let $E \in \real^{n_x \times |\mathcal{E}|}$ be the  oriented node-edge incidence matrix \cite[Chapter 8]{FB:16} of a weakly connected acyclic\footnote{The acyclic assumption implies that $\mathrm{rank}(E^{\sf T}) = |\mathcal{E}|$, in line with our assumption from Section \ref{Sec:ReviewPrimal} that the constraint matrix has full row rank. This assumption can be relaxed at the expense of more complex stability/performance proofs.} graph $\mathcal{G} = (\{1,\ldots,n_x\},\mathcal{E})$, where $\mathcal{E}$ is the set of oriented edges. The dual problem $\maximize_{\nu \in \real^{n_r}}\Phi(\nu)$ is then equivalent to the constrained problem\footnote{Another equivalent formulation can be obtained by using the Laplacian matrix $EE^{\sf T} = {\rm L} = {\rm L}^{\sf T} \in \real^{n_x \times n_x}$ of the graph, and using instead the constraint $({\rm L} \otimes I_{n_r})\boldsymbol{\nu} = \vzeros[{(n_x n_r)}]$. This formulation is sometimes preferable for multi-agent implementations, the analysis of which requires only small modifications from the present analysis. We focus instead of formulations involving the incidence matrix.}
\begin{equation}\label{Eq:DualEqualityConstrained}
\begin{aligned}
& \underset{\boldsymbol{\nu} \in {\mathbb{R}^{(n_x n_r)}}}{\text{minimize}}
& & J_{\rm dual}(\boldsymbol{\nu}) := -\sum_{i=1}^{n_x} \widetilde{\Phi}_i(\nu^i) \\
& \text{subject to}
&& (E^{\sf T} \otimes I_{n_r})\boldsymbol{\nu} = {\vzeros[|\mathcal{E}|* n_r]}\,,
\end{aligned}
\end{equation}
where $\boldsymbol{\nu} = \mathrm{col}({\nu^1},\nu^2,\ldots,{{\nu^{n_x}}}) \in \real^{{(n_x n_r)}}$. Since the graph is acyclic, $E^{\sf T}$ has full row rank, and therefore satisfies our assumption concerning the constraint matrix (Section \ref{Sec:PrimalDualMethods}). The key observation now is that the parameters $b$ do not enter into the equality constraints of the optimization problem \eqref{Eq:DualEqualityConstrained}; this permits the full application of augmented Lagrangian techniques for improving the $\mathcal{H}_2$ performance of the saddle point algorithm. Building the augmented Lagrangian \eqref{Eq:LAug} for the problem \eqref{Eq:DualEqualityConstrained} we have
$$
L_{\rm dual}^{\rm aug} = -\sum_{i=1}^{n_x} \widetilde{\Phi}_i(\nu^i) + \boldsymbol{\mu}^{\sf T}(E^{\sf T} \otimes I_{n_r})\boldsymbol{\nu} + \frac{\rho}{2}\boldsymbol{\nu}^{\sf T}({\rm L} \otimes I_{n_r})\boldsymbol{\nu}\,,
$$
where $\boldsymbol{\mu} = \mathrm{col}({\mu^1},\mu^2,\ldots,{\mu^{|\mathcal{E}|}}) \in {\real^{|\mathcal{E}|* n_r}}$ is a stacked vector of Lagrange multipliers $\mu^\ell \in \real^{n_r}$ for $\ell \in \{1,\ldots,|\mathcal{E}|\}$, and ${\rm L} = {\rm L}^{\sf T} = EE^{\sf T} \in \real^{n_x \times n_x}$ is the Laplacian matrix of the undirected graph $\mathcal{G}_{\rm u}$, obtained by ignoring the orientation of the edges in $\mathcal{G}$; we let $\mathcal{N}(i)$ denote the neighbours of vertex $i$ in the graph $\mathcal{G}_{\rm u}$. Applying the saddle-point method, the dynamics may be written block-component-wise as
\begin{equation}\label{Eq:DecentralizedDual}
\begin{aligned}
\bar{\tau}_{\nu}\dot{\nu}^{i} &= -\frac{s_is_i^{\sf T}}{q_i}\nu^i - \left(\frac{c_i}{q_i}s_i+w_ib_i\right)\\
&\quad - \sum_{j:(i,j)\in\mathcal{E}}\mu^{ij}+\sum_{j:(j,i)\in\mathcal{E}}\mu^{ji}\\
&\quad - \rho \sum_{j \in \mathcal{N}(i)}\nolimits (\nu^i-\nu^j)\,, \quad i \in \{1,\ldots,n_x\}\,,\\
\bar{\tau}_{\mu}\dot{\mu}^{ij} &= \nu^i - \nu^j\,, \quad (i,j) \in \mathcal{E}\,,\\
\end{aligned}
\end{equation}
The algorithm \eqref{Eq:DecentralizedDual} is distributed, in that the $i$th update equation requires only the local parameters $s_i, w_i, b_i, c_i, q_i$ along with communicated state variables $\nu^j, \mu^{ij}$ which come from adjacent nodes and edges in the graph $\mathcal{G}$. We refer to this dynamical system as the augmented dual distributed saddle-point (ADD-SP) dynamics. Following \eqref{Eq:DualDynamics}, we equip this system with disturbance inputs and performance outputs
\begin{equation}\label{Eq:ADD-SP-Performance-Outputs}
b_i {\mapsto} b_i + t_b \eta_i\,, \quad z_i = -q_i^{-\frac{1}{2}}s_i^{\sf T}\nu^i\,, \quad i \in \{1,\ldots,n_x\}\,.
\end{equation}
While a closed-form expression for the $\mathcal{H}_2$ norm of the this system is difficult to compute, we can state the following comparative result.

\medskip

\begin{corollary}[\bf Augmented Dual Distributed Saddle-Point Performance]\label{Cor:ADD}
Consider the ADD-SP dynamics \eqref{Eq:DecentralizedDual}, denoted by $G_{\rm ADD}$, with disturbance inputs $\eta(t)$ and performance output $z(t)$ as in \eqref{Eq:ADD-SP-Performance-Outputs}, under the same assumptions as Theorem \ref{Thm:AugmentedPrimalOutputs}. Then the squared $\mathcal{H}_2$ norm of the system \eqref{Eq:DecentralizedDual}-\eqref{Eq:ADD-SP-Performance-Outputs} satisfies the upper bound
\begin{equation}\label{Eq:GAddNorm}
\|G_{\rm ADD}\|_{\mathcal{H}_2}^2 \leq \frac{t_b^2}{2\bar{\tau}_{\nu}}\mathrm{Tr}(\mathcal{W}_b^{\sf T}\mathcal{W}_b)\,,
\end{equation}
where $\mathcal{W}_b = \mathrm{blkdiag}(w_1,w_2, \ldots, w_{n_x}) \in \real^{(n_x n_r) \times n_x}$.
Moreover, \eqref{Eq:GAddNorm} is satisfied with equality if and only if $\rho = 0$.
\end{corollary}

\begin{pfof}{Corollary \ref{Cor:ADD}}
To begin, we note that
$$
\sum_{i=1}^{n_x} \widetilde{\Phi}_i(\nu^i) = -\frac{1}{2}\boldsymbol{\nu}^{\sf T}\mathcal{S}Q^{-1}\mathcal{S}^{\sf T}\boldsymbol{\nu} - \boldsymbol{\nu}^{\sf T}\mathcal{S}Q^{-1}c - {\boldsymbol{\nu}^{\sf T}\mathcal{W}_b b}
$$
where $\mathcal{S} = \mathrm{blkdiag}(s_1,s_2,\ldots,s_{n_x}) \in \real^{(n_x n_r) \times n_x}$.

In vector notation after shifting the equilibrium to the origin, the system \eqref{Eq:DecentralizedDual}--\eqref{Eq:ADD-SP-Performance-Outputs} can be written as
$$
\begin{aligned}
\bar{\tau}_{\nu}\dot{\boldsymbol{\nu}} &= -\mathcal{S}Q^{-1}\mathcal{S}^{\sf T}\boldsymbol{\nu} - t_b\mathcal{W}_b\eta_b - (E \otimes I_{n_r})\boldsymbol{\mu} - \rho ({\rm L}\otimes I_{n_r})\boldsymbol{\nu}\\
\bar{\tau}_{\mu}\dot{\boldsymbol{\mu}} &= (E^{\sf T} \otimes I_{n_r})\boldsymbol{\nu}\\
z &= -Q^{-\frac{1}{2}}\mathcal{S}^{\sf T}\boldsymbol{\nu}
\end{aligned}
$$
where $\eta = \mathrm{col}(\eta_1,\ldots,\eta_{n_x}) \in \real^{n_x}$. In state-space this translates to
\begin{equation}\label{Eq:RA4_gen}
\begin{aligned}
	\begin{bmatrix}
		\dot{\boldsymbol{\nu}}\\
		\dot{\boldsymbol{\mu}}
	\end{bmatrix}
&=
	\underbrace{\begin{bmatrix}
		-\frac{1}{\bar{\tau}_{\nu}}\mathcal{S}Q^{-1}\mathcal{S}^{\sf T}- \frac{1}{\bar{\tau}_{\nu}}\rho ({\rm L}\otimes I_{n_r})& -\frac{1}{\bar{\tau}_{\nu}} (E \otimes I_{n_r})\\
		\frac{1}{\bar\tau_\mu} (E^{\sf T} \otimes I_{n_r})& \vzeros[]
	\end{bmatrix}}_{A_\text{ADD}}
	\begin{bmatrix}
		\boldsymbol{\nu}\\
		\boldsymbol{\mu}
	\end{bmatrix}\\
	& -
	\underbrace{\begin{bmatrix}\frac{t_b}{\bar\tau_\nu}\mathcal{W}_b\\\vzeros[]\end{bmatrix}}_{B_\text{ADD}}
	\eta_b, \, \, \quad
z = \underbrace{\begin{bmatrix}-Q^{-\frac{1}{2}}\mathcal{S}^{\sf T} & \vzeros[] \end{bmatrix}}_{C_\text{ADD}}\begin{bmatrix}\boldsymbol{\nu} \\ \boldsymbol\mu \end{bmatrix}\,.
\end{aligned}
\end{equation}
As in our previous results, one may verify that $A_{\rm ADD}$ is Hurwitz and that $(A_{\rm ADD},C_{\rm ADD})$ is observable. Consider the observability Gramian candidate $X_\text{ADD} = \frac{1}{2}\mathrm{blkdiag}(\bar\tau_\nu\, I_{(n_x n_r)},\bar\tau_{\mu}\, {I_{|\mathcal{E}| * n_r}})$. Straightforward algebra shows that
\begin{equation}
\label{Eq:regGramian_RA4}
X_\text{ADD}\, A_\text{ADD}+A_\text{ADD}^{\sf T}X_\text{ADD}+Q_\text{ADD}= \begin{bmatrix}-\rho ({\rm L}\otimes I_{n_r}) & \vzeros[]\\ \vzeros[] & \vzeros[]\end{bmatrix} \preceq \vzeros [],
\end{equation}
where $Q_\text{ADD}=C_\text{ADD}^\top C_\text{ADD}$.
We conclude that $X_\text{ADD}$ is a generalized observability Gramian for the system $G_\text{ADD}$. Furthermore, if $X_{t}$ is the true observability Gramian for $G_{\rm ADD}$, then as in Theorem~\ref{Thm:PrimalOutputs_reg}, 
we have $X_{t} \preceq X_\text{ADD}$, and we conclude that
\begin{equation}\label{Eq:RA4_rho}
\|G_\text{ADD}\|_{\mathcal{H}_2}^2 \leq \mathrm{Tr}(B_\text{ADD}^{\sf T}X_\text{ADD}B_\text{ADD})=\frac{t_b^2}{2 \bar\tau_\nu}\, \mathrm{Tr}(\mathcal{W}_b^{\sf T}\mathcal{W}_b)\,.%
\end{equation}
When $\rho = 0$, $X_{\rm ADD}$ is the exact observability Gramian and hence 
\begin{equation}\label{Eq:NimaEquality}
\|G_\text{ADD}\|_{\mathcal{H}_2}^2=\frac{t_b^2}{2 \bar\tau_\nu}\mathrm{Tr}(\mathcal{W}_b^{\sf T}\mathcal{W}_b)\,
\end{equation}

To complete the proof, it remains to show that \eqref{Eq:NimaEquality} in fact implies that $\rho=0$. Suppose that \eqref{Eq:NimaEquality} holds, and define $\Delta := X_\text{ADD}-X_t \succeq \vzeros[]$. Then, following a similar  argument as in the proof of  Theorem \ref{Thm:PrimalOutputs_reg}, one may show that $\Delta B_\text{ADD} = \vzeros[]$. 
From \eqref{Eq:regGramian_RA4} and the fact that $X_t$ is the actual observability Gramian of the system, we may subtract equations to obtain
\[
\Delta\, A_\text{ADD}+A_\text{ADD}^{\sf T}\Delta= \begin{bmatrix}-\rho ({\rm L}\otimes I_{n_r}) & \vzeros[]\\ \vzeros[] & \vzeros[]\end{bmatrix}\,. \]
Now, by multiplying each side of the equality above from the left by $B_\text{ADD}^{\sf T}$ and from the right by $B_\text{ADD}$, we find that 
$$
\begin{aligned}
\vzeros[] &= \rho \mathcal{W}_{b}^{\sf T}({\rm L}\otimes I_{n_r})\mathcal{W}_b\\
&= \rho \mathcal{W}_{b}^{\sf T}(E \otimes I_{n_r})(E^{\sf T} \otimes I_{n_r})\mathcal{W}_b = \rho M^{\sf T}M
\end{aligned}
$$
where $M = (E^{\sf T} \otimes I_{n_r})\mathcal{W}_b$. Since the graph $\mathcal{G}$ is connected and $W_b$ is square and of full rank, it always holds that $M \neq \vzeros[]$, and therefore we conclude that 
$\rho=0$.
\end{pfof}

The point of interest from Corollary \ref{Cor:ADD} is that the bound on the $\mathcal{H}_2$ performance is \emph{independent} of $\rho$. In particular then, as $\rho$ becomes large the norm does not grow without bound, which resolves the issue observed in the result of Theorem \ref{Thm:AugmentedPrimalOutputs}. When applied to the resource allocation problem in Section \ref{Sec:Resource}, we will in fact be able to make the stronger statement that the norm is a strictly decreasing function of $\rho$, and augmentation can therefore be used successfully to improve saddle-point algorithm performance. 

\section{Application to Resource Allocation Problems}
\label{Sec:Resource}

We now apply the results from the previous subsections to a particular class of problems. As a special case of the problem \eqref{Eq:ConvexEqualityConstrained}, consider the \emph{resource allocation problem}
\begin{equation}\label{Eq:Resource}
\begin{aligned}
& \underset{x \in \mathbb{R}^n}{\text{minimize}}
& & \sum_{i=1}^{n}\nolimits \frac{1}{2}q_ix_i^2 + c_ix_i \\
& \text{subject to}
&& \sum_{i=1}^{n}\nolimits x_i = \sum_{i=1}^n\nolimits d_i\,,
\end{aligned}
\end{equation}
where $q_i > 0$, $c_i \in \real$, and $d_i \in \real$. Comparing \eqref{Eq:Resource} to \eqref{Eq:ConvexEqualityConstrained}, we have $Q = \mathrm{diag}(q_1,\ldots,q_n)$, $S = \vones[n]^{\sf T}$, $W_b = \vones[n]^{\sf T}$, and $d := \mathrm{col}(d_1,\ldots,d_n) =  b$. The interpretation of \eqref{Eq:Resource} is that a resource must be obtained from one of $n$ suppliers in amount $x_i$, subject to a total demand satisfaction constraint. The objective function of  \eqref{Eq:Resource} can be interpreted as the sum of the utilities $-c_ix_i$ minus the sum of the costs $q_ix_i^2/2$.
In a multi-agent context, each variable $x_i$ is assigned to an agent, the parameters $q_i, c_i, d_i$ are available locally to each agent, and the agents must collectively solve the problem \eqref{Eq:Resource} through local exchange of information. 
As a concrete example, in the context of power system frequency control, the objective function models the cost of producing an auxiliary power input $x_i$; see our preliminary work \cite{JWSP-BKP-NM-FD:16c} for additional details.

We will consider \eqref{Eq:Resource} along with several equivalent reformulations, and apply our results from the previous sections to assess the input-output performance of the resulting saddle-point algorithms. External disturbances $\eta_{d}$ will be integrated into the algorithms as $d_i \mapsto d_i + \eta_i$, where $\eta_i$ models the disturbances in demand $d_i$. For simplicity, all time-constant matrices are assumed to be multiples of the identity matrix. To most clearly indicate which algorithms require communication of which variables, in this section algorithms are \emph{not} written in deviation coordinates with respect to the optimizer. In all cases, the performance output $z$ is chosen such that $\|z(t)\|_2^2 = (x(t)-x^{\star})^{\sf T}Q(x(t)-x^{\star})$, where $x^{\star}$ is the global primal optimizer. With this choice of performance output, the transient performance of the algorithm is measured using the same relative weightings as the steady-state performance.

To begin, the augmented Lagrangian of \eqref{Eq:Resource} is given in vector notation by
\begin{equation}\label{Eq:RALagrangian}
L(x,\nu) = \frac{1}{2}x^{\sf T}Qx + c^{\sf T}x + \nu \vones[n]^{\sf T}(x-d) + \frac{\rho}{2}\|\vones[n]^{\sf T}(x-d)\|_2^2\,.
\end{equation}
where $\nu \in \real$. Applying the saddle-point method to the Lagrangian $L(x,\nu)$ and attaching the same disturbances $\eta \in \real^n$ and performance outputs $z \in \real^n$ as before, we obtain the {\em centralized saddle-point} dynamics
\begin{equation}\label{Eq:ResAlccent}
\mathrm{RA}_{\rm cent}(\rho): 
\begin{cases}
\begin{aligned}
\bar{\tau}_x \dot{x} &= -Qx - c - \nu \vones[n] - \rho \vones[n]\vones[n]^{\sf T}(x-d-\eta)\\
\bar{\tau}_{\nu}\dot{\nu} &= \vones[n]^{\sf T}(x - d - \eta)\\
z &= Q^{\frac{1}{2}}(x-x^*).
\end{aligned}
\end{cases}
\end{equation}

When $\rho = 0$, the algorithm \eqref{Eq:ResAlccent} is of a gather-and-broadcast type \cite{FS-SG:17}, where all states $x_i$ and disturbances $d_i$ are collected and processed by a central agent with state $\nu$. When $\rho > 0$, the additional term $\vones[n]\vones[n]^{\sf T}$ in the algorithm requires all-to-all communication of the local imbalances $x_i - d_i$. 

{We now consider a reformulation that results in a distributed optimization algorithm.} Let $\mathcal{G} = (\{1,\ldots,n\},\mathcal{E})$ denote a weakly connected acyclic graph over the agent set $\{1,\ldots,n\}$, and let $E \in \real^{n \times |\mathcal{E}|}$ denote the oriented node-edge incidence matrix of $\mathcal{G}$. The constraint $\vones[n]^{\sf T}x= \vones[n]^{\sf T}d$ in the resource allocation problem \eqref{Eq:Resource} is equivalent to the existence\footnote{This follows since $\mathrm{ker}(E^{\sf T}) = \mathrm{span}(\vones[n])$, and hence $\mathrm{Im}(E)$ is the subspace orthogonal to the vector $\vones[n]$.} of a vector $\delta \in \real^{|\mathcal{E}|}$ such that $E\delta = x - d$.
The resource allocation problem \eqref{Eq:Resource} can therefore be equivalently rewritten as
\begin{equation}\label{Eq:ResourceDist}
\begin{aligned}
& \underset{x \in \real^n, \delta \in \real^{|\mathcal{E}|}}{\text{minimize}}
& & \sum_{i=1}^{n}\nolimits \frac{1}{2}q_ix_i^2 + c_ix_i \\
& \text{subject to}
&& E\delta = x - d\,,
\end{aligned}
\end{equation}
with associated augmented Lagrangian
$$
L^\prime(x,\delta,\nu) = \frac{1}{2}x^{\sf T}Qx + c^{\sf T}x + \nu^{\sf T}(E\delta - x + d) + \frac{\rho}{2}\|E\delta - x + d\|_2^2\,,
$$
where $\nu \in \real^n$. This reformulation can be interpreted as a version of \eqref{Eq:ConvexEqualityConstrained} with an expanded primal variable $(x,\delta)$ and an expanded dual variable $\nu = \mathrm{col}(\nu_1,\ldots,\nu_n)$. 
By applying the saddle-point method to the problem \eqref{Eq:ResourceDist}, we obtain the \emph{distributed saddle-point} dynamics
\begin{equation}\label{Eq:ResAlcdist}
\mathrm{RA}_{\rm dist}(\rho): 
\begin{cases}
\begin{aligned}
\bar{\tau}_x \dot{x} &= -Qx - c + \rho(E\delta-x+d + \eta) + \nu\\
\bar{\tau}_{\delta}\dot{\delta} &= -E^{\sf T}\nu - \rho E^{\sf T}(E\delta-x+d+\eta)\\
\bar{\tau}_{\nu}\dot{\nu} &= E\delta - x + d + \eta\\
z &= Q^{\frac{1}{2}}(x-x^*).
\end{aligned}\,
\end{cases}
\end{equation}
When $\rho = 0$, the algorithm \eqref{Eq:ResAlcdist} is distributed with the topology of the graph $\mathcal{G}$, with states $(x_i,\nu_i)$ associated with each node and a state $\delta_i$ associated with each edge. When $\rho > 0$, the algorithm contains the so-called edge Laplacian matrix $E^{\sf T}E$ \cite{DZ-MM:11}, which under our acyclic assumption is positive definite.

Our third formulation is the dual ascent algorithm \eqref{Eq:DualFunction}. Substitution of the appropriate matrices into \eqref{Eq:DualFunction} leads to the \emph{centralized dual ascent} dynamics
\begin{equation}\label{Eq:ResAlcdual}
\mathrm{RA}_{\rm cent}^{\rm dual}: 
\begin{cases}
\begin{aligned}
\bar{\tau}_{\nu}\dot{\nu} &= -(\vones[n]^{\sf T}Q^{-1}\vones[n])\nu - \vones[n]^{\sf T}(Q^{-1}c+d+\eta)\\
z &= Q^{\frac{1}{2}}(x-x^*) = -Q^{-\frac{1}{2}}(c + \nu \vones[n]) - Q^{\frac{1}{2}}x^*\,,
\end{aligned}
\end{cases}
\end{equation}
where $\nu \in \real$. Algorithm \eqref{Eq:ResAlcdual} is again centralized, with a single central agent with state $\nu$ performing all computations and broadcasting $x_i = -q_i^{-1}(c_i+\nu)$ back to each agent.
For our fourth and final formulation, we apply the ADD-SP method developed in Section \ref{Sec:DistributedDual}. For the problem \eqref{Eq:Resource}, one quickly deduces that $\mathcal{S} = \mathcal{W}_b = I_n$, and the algorithm \eqref{Eq:DecentralizedDual} reduces to
\begin{equation}\label{Eq:ResAlcdist1}
\mathrm{RA}_{\rm dist}^{\rm dual}(\rho): 
\begin{cases}
\begin{aligned}
\bar{\tau}_\nu \dot{\nu} &= -Q^{-1}\nu - (d+\eta) - Q^{-1}c - E\mu - \rho {\rm L} \nu\\
\bar{\tau}_{\mu}\dot{\mu} &= E^{\sf T}\nu\\
z &= Q^{\frac{1}{2}}(x-x^*) = -Q^{-\frac{1}{2}}(\nu+c)+Q^{\frac{1}{2}}x^*\,,
\end{aligned}
\end{cases}
\end{equation}
When $\rho = 0$, this algorithm is distributed with the graph $\mathcal{G}$ associated with the incidence matrix $E$, with states $\nu_i$ associated with nodes and states $\mu_{ij}$ associated with edges. When $\rho > 0$, the algorithm additionally contains the undirected Laplacian matrix ${\rm L} = EE^{\sf T}$ of $\mathcal{G}$, and thus remains distributed.

For each of the four formulations above, we compute the $\mathcal{H}_2$ norm from the disturbance input $\eta$ to the performance output $z$.
For $\mathrm{RA}_{\rm cent}(\rho), \mathrm{RA}^{\rm dual}_{\rm cent}$, and $\mathrm{RA}^{\rm dual}_{\rm dist}(\rho)$ this follows immediately from Theorem \ref{Thm:AugmentedPrimalOutputs}, Proposition \ref{Prop:Dual-Ascent}, and Corollary \ref{Cor:ADD}, respectively. The algorithm $\mathrm{RA}_{\rm dist}(\rho)$ requires a modification of the proof of Theorem \ref{Thm:AugmentedPrimalOutputs}, since the objective function is no longer strongly convex in the primal variables $(x,\delta)$; we omit the details.

\begin{table}[h!]
  \centering
  \caption{Comparison of squared $\mathcal{H}_2$ norm expressions}
  \label{Tab:1}
  \begin{tabular}{l|cc}
   \hline
    System & $\rho = 0$	& $\rho \rightarrow \infty$\\ [0.5ex]
    \hline
    $\mathrm{RA}_{\rm cent}(\rho)$ & $n/(2\bar{\tau}_{\nu})$   & $+\infty$\\
  \\
    $\mathrm{RA}_{\rm dist}(\rho)$ & $n/(2\bar{\tau}_{\nu})$ & $+\infty$\\
        \\
    $\mathrm{RA}_{\rm cent}^{\rm dual}$ & $n/(2\bar{\tau}_{\nu})$ & independent of $\rho$\\
        \\
    $\mathrm{RA}_{\rm dist}^{\rm dual}(\rho)$ & $n/(2\bar{\tau}_{\nu})$ & $1/(2\bar{\tau}_{\nu})$\\ [0.5ex]
    \hline
  \end{tabular}
\end{table}

\begin{figure*}
\centering
\includegraphics[height=1\columnwidth]{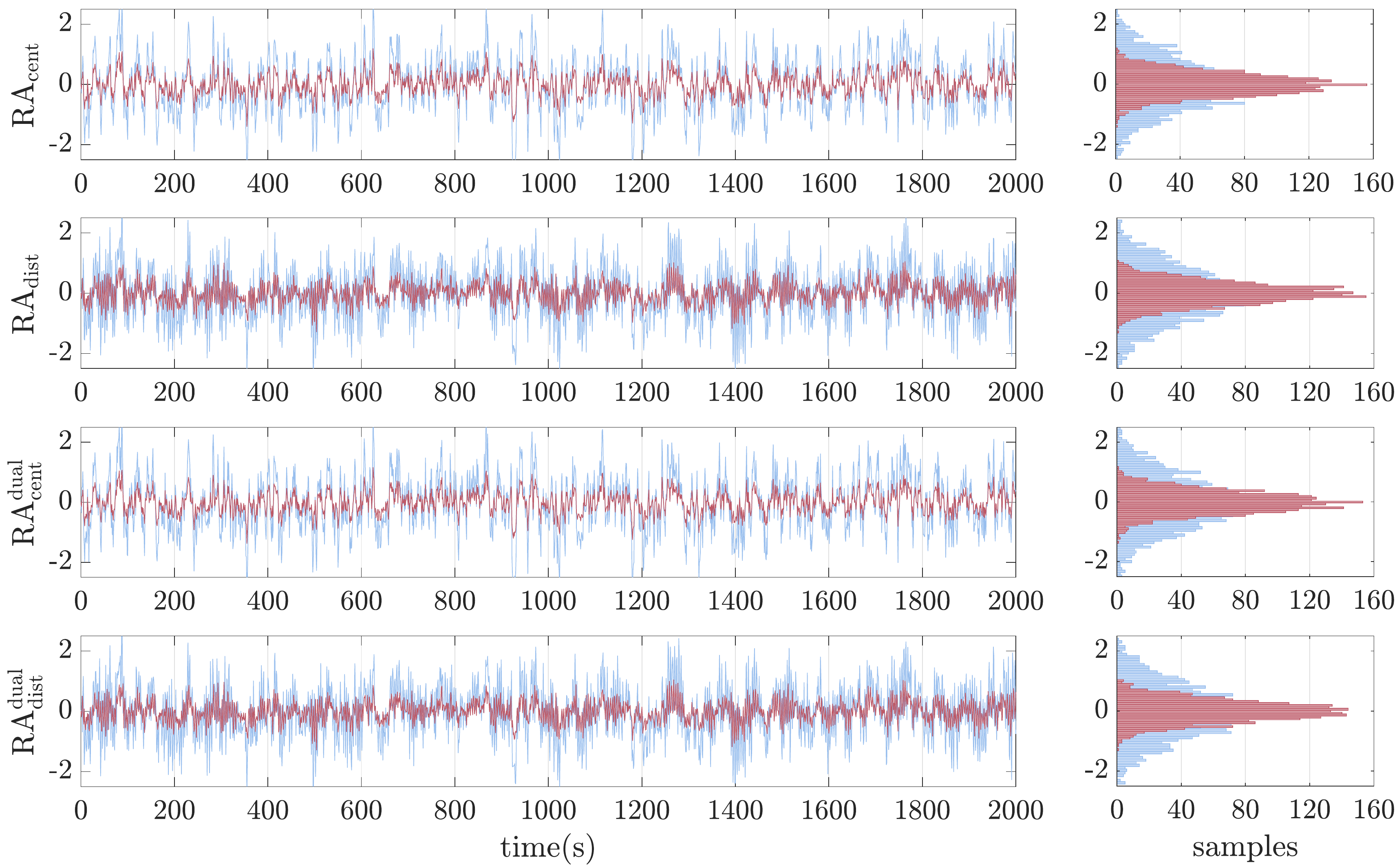}
\caption{Steady-state variance for the four different implementations for the unaugmented case, for parameters $n=2$, $Q=\diag(4,25)$, $\tau_x=\tau_\delta=\tau_\nu=\tau_\mu=1$, $E=[1\,\, -1]^\top$; the remaining parameters do not influence the results.}
\label{Fig:UnAug}
\end{figure*}

\begin{figure*}
\centering
\includegraphics[height=0.765\columnwidth]{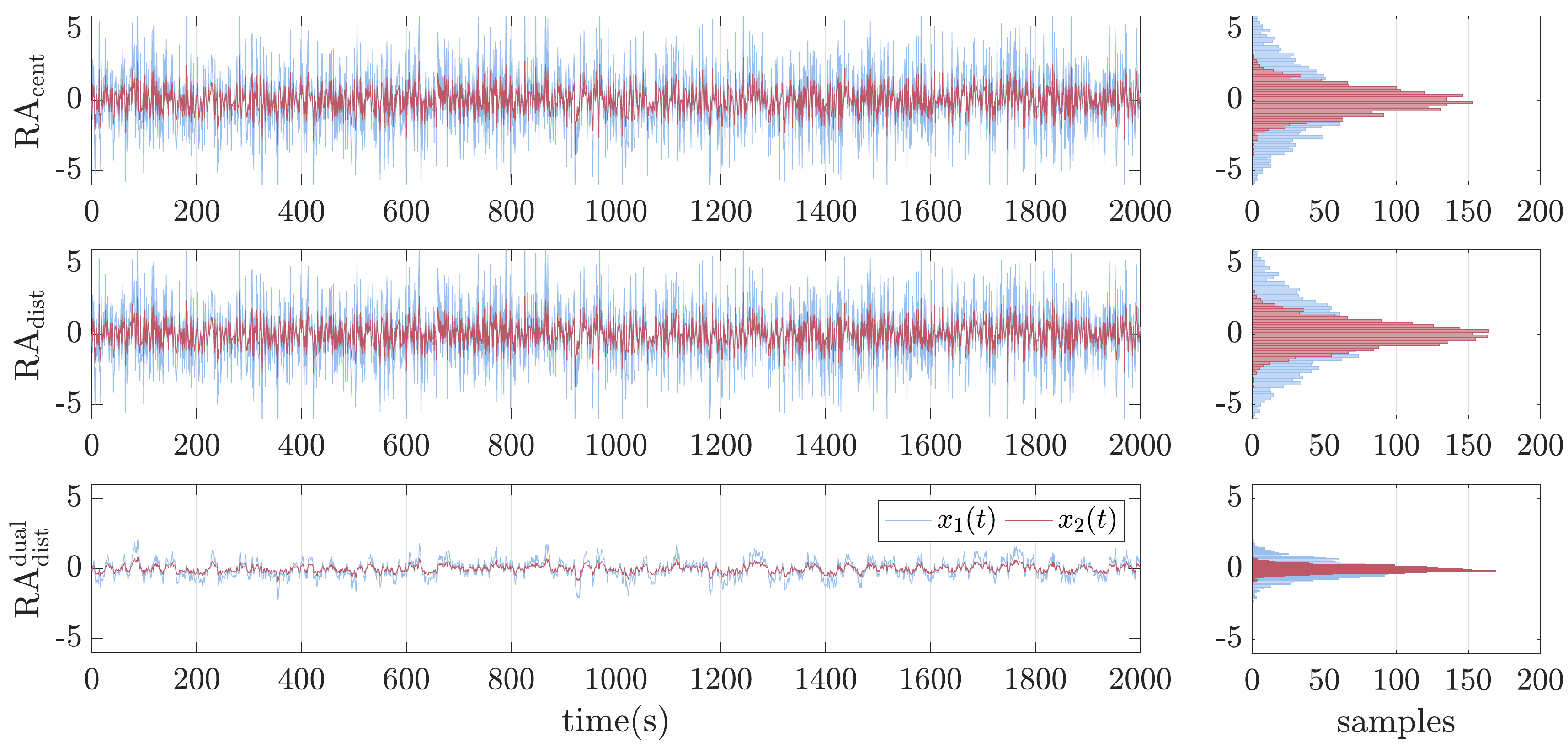}
\caption{Steady-state variance for the three different implementations for the augmented case, for parameters $n=2$, $\rho=100$, $Q=\diag(4,25)$, $\tau_x=\tau_\delta=\tau_\nu=\tau_\mu=1$, $E=[1 -1]^\top$; the remaining parameters do not influence the results.}
\label{Fig:Aug}
\end{figure*}

The first column of Table \ref{Tab:1} shows the $\mathcal{H}_2$ system norms for the four formulations when $\rho = 0$, i.e., the unaugmented versions of the various saddle-point algorithms. Despite substantial differences between the algorithms in terms of information structure and number of states, all four have the same input-output performance in the $\mathcal{H}_2$ norm. This implies that a distributed implementation will perform no worse than a centralized implementation. 

\begin{figure*}[!ht]
        \begin{center}
        \begin{subfigure}[!ht]{\columnwidth}
        		\centering
                \includegraphics[height=0.85\columnwidth]{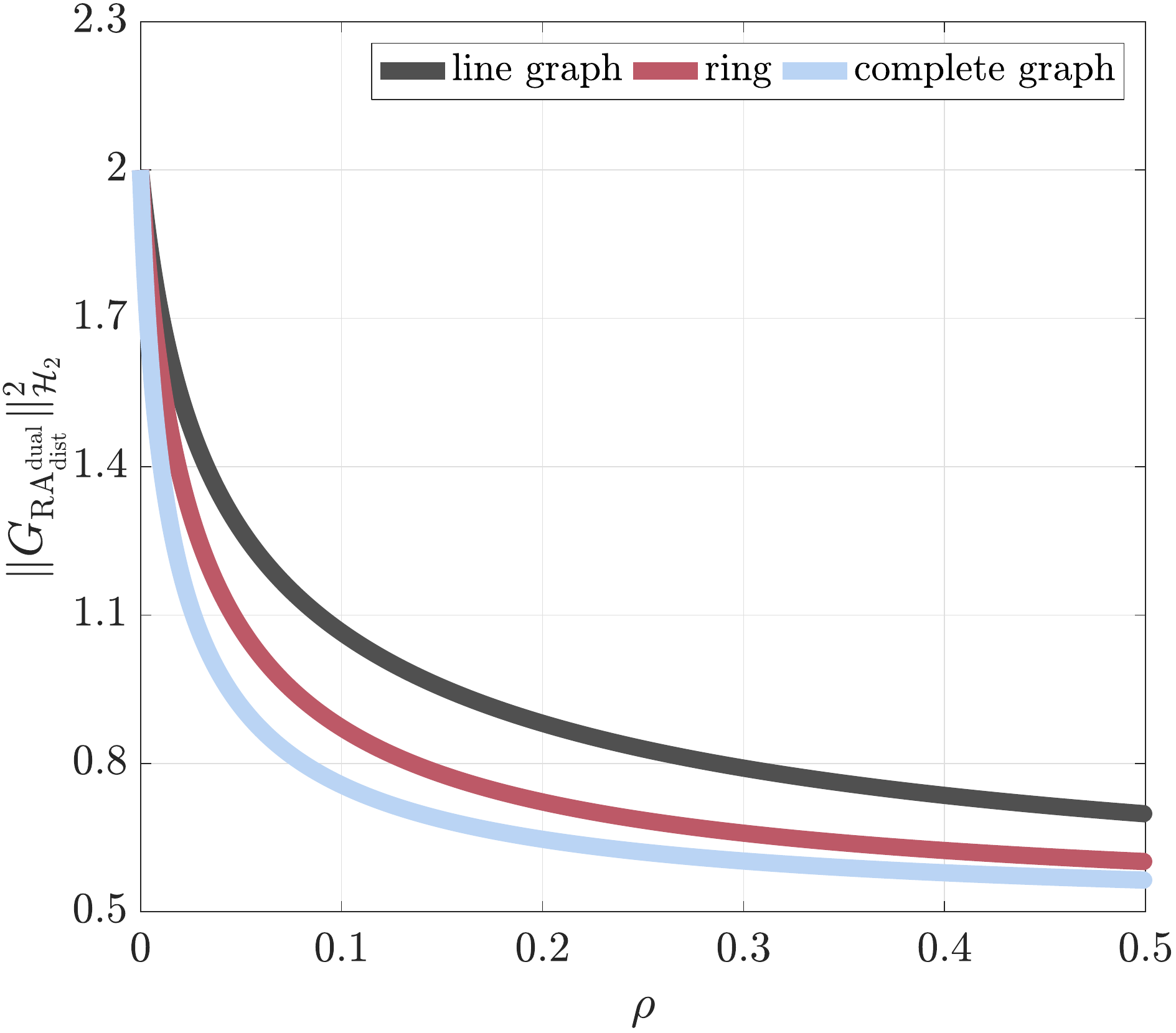}
                \caption[]{$\mathrm{RA}_{\rm dist}^{\rm dual}(\rho)$ for $Q=\diag(4,25,16,49)$}
                \label{Fig:RA4}
        \end{subfigure}%
        \begin{subfigure}[!ht]{\columnwidth}
        		\centering
                \includegraphics[height=0.85\columnwidth]{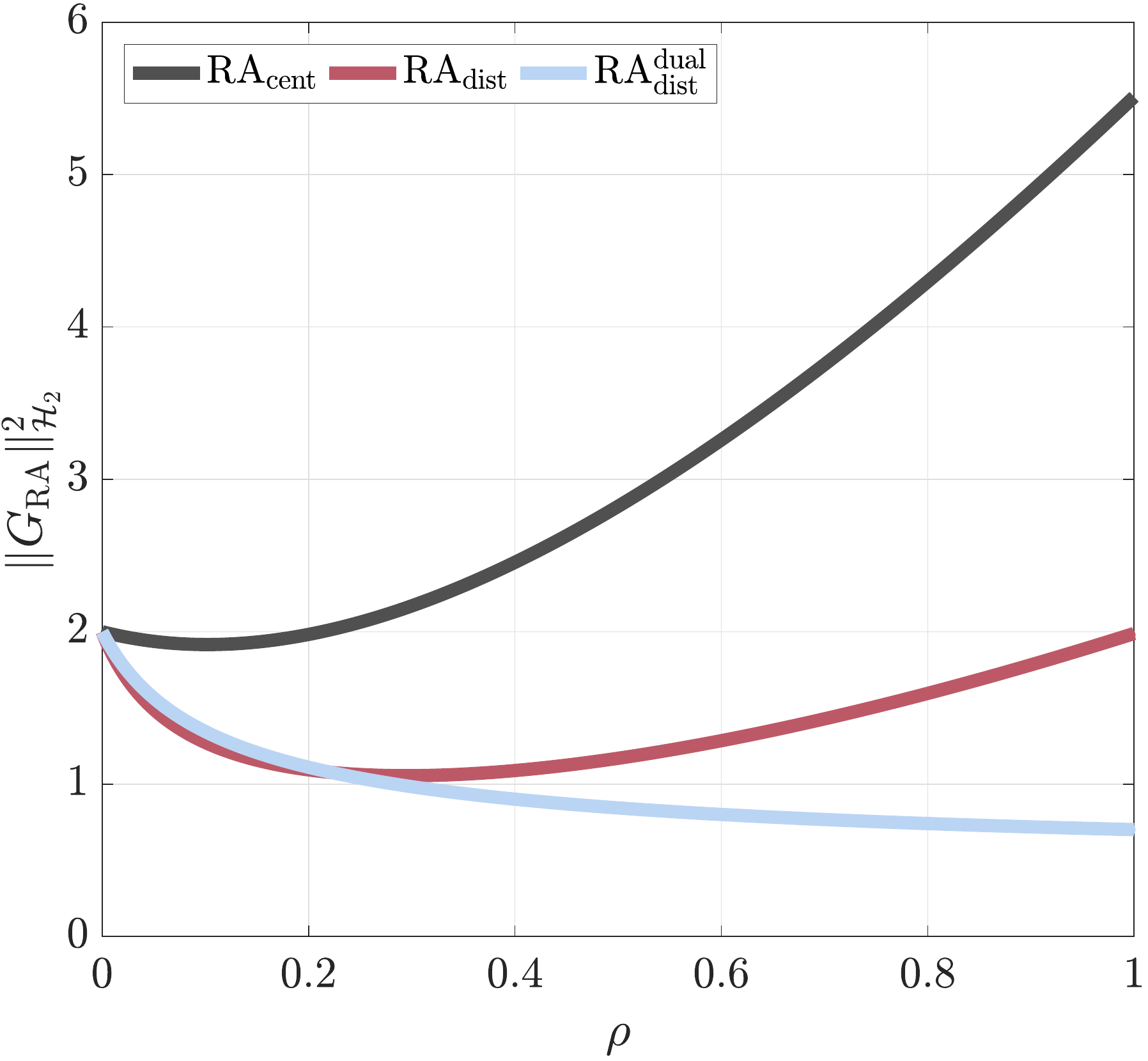}
                 \caption[]{$\mathrm{RA}_{\rm cent}$, $\mathrm{RA}_{\rm dist}$, $\mathrm{RA}_{\rm dist}^{\rm dual}$ for $Q=\diag(4,4,4,9)$ and  line graph}
                \label{Fig:RA}
                \end{subfigure}
		\caption[]{Squared $\mathcal{H}_2$ norm as a function of $\rho$, for parameters $n=4$, $\tau_x=\tau_\delta=\tau_\nu=\tau_\mu=1$ and unweighted graphs.}
		\label{Fig:AugImp}
\end{center}
\end{figure*}

While these four formulations all possess identical system norms under the basic primal-dual algorithm, augmentation differentiates these methods from one another, and substantial differences between the algorithms begin to appear as $\rho$ is increased. The limiting results are tabulated in the second column of Table \ref{Tab:1}. The input-output performance of the first two formulations becomes arbitrarily bad as the augmentation gain $\rho$ increases, while the performance of the ADD-SP algorithm improves substantially, becoming independent of the system size in the limit $\rho \to \infty$. 

We illustrate the results in Table \ref{Tab:1} via time-domain simulations in Figures~\ref{Fig:UnAug} and \ref{Fig:Aug} for the system in \eqref{Eq:ResourceDist} with $n=2$ and an underlying line graph with $E=[1\, -1]^\top$. With unit variance white noise as inputs, the unaugmented implementation in Figure~\ref{Fig:UnAug} for the four different algorithms results in identical steady-state output variance, numerically computed as the squared $\mathcal{H}_2$ norm in \eqref{Eq:H2norm}. 

In Figure~\ref{Fig:Aug}, a sufficiently large augmentation factor $\rho$ is introduced to penalise the constraint violations. It is observed that with the augmentation, the steady-state variance of the outputs for the centralized and distributed implementations in $\mathrm{RA}_{\rm cent}(\rho)$, $\mathrm{RA}_{\rm dist}(\rho)$ worsens, while that of the distributed dual implementation from $\mathrm{RA}_{\rm dist}^{\rm dual}(\rho)$ improves. 

Figure \ref{Fig:AugImp}(a) illustrates how the choice of communication graph topology influences the performance of the algorithm $\mathrm{RA}_{\rm dist}^{\rm dual}(\rho)$. We consider $n = 4$ agents, and implement the algorithm with line, ring, and complete communication graphs. While Corollary \ref{Cor:ADD} and the results of Table \ref{Tab:1} hold only for acyclic graphs, Figure \ref{Fig:AugImp}(a) shows that in all cases the algorithm's performance improves as $\rho$ increases. For a given value of $\rho$, graphs with higher connectivity show a greater improvement. This behaviour is explained by noting that the algorithm \eqref{Eq:ResAlcdist1} has the same form as the augmented saddle-point dynamics \eqref{Eq:SaddlePointAugmented}, and an analysis similar to that performed for Theorem \ref{Thm:AugmentedPrimalOutputs} can in fact be performed for \eqref{Eq:ResAlcdist1}. For uniform cost function parameters and an acyclic graph, this leads to the expression
\begin{align}\label{Eq:Whatever}
\|G_{\mathrm{RA}_{\rm dist}^{\rm dual}(\rho)}\|_{\mathcal{H}_2}^2 &= \frac{1}{2\bar{\tau}_{\nu}}\left(1 + \sum_{i=2}^{n}\frac{q}{q+\rho \lambda_i}\right),
\end{align}
where $0 = \lambda_1 < \lambda_2 \leq \cdots \leq \lambda_n$ are the eigenvalues of the Laplacian matrix. As graph connectivity increases, so does $\lambda_2$, and performance therefore improves. From a design perspective, note that for a fixed time constant $\bar{\tau}_{\nu}$ and a desired level of performance $\gamma \in (1/\sqrt{2\bar{\tau}_{\nu}},\sqrt{n}/\sqrt{2\bar{\tau}_{\nu}})$, examination of \eqref{Eq:Whatever} shows that a sufficient condition for $\|G_{\mathrm{RA}_{\rm dist}^{\rm dual}(\rho)}\|_{\mathcal{H}_2} \leq \gamma$ is that
$$
\rho \geq \frac{q}{\lambda_2}\frac{n-2\bar{\tau}_{\nu}\gamma^2}{2\bar{\tau}_{\nu}\gamma^2-1}\,.
$$
In particular, this shows that the augmentation gain should be chosen in inverse proportion to the algebraic connectivity $\lambda_2$ of the Laplacian matrix. A strongly connected graph will therefore require a lower augmentation gain than a weakly connected graph to achieve a desired level of $\mathcal{H}_2$ performance. Achieving an $\mathcal{H}_{2}$ norm lower than $1/\sqrt{2\bar{\tau}_{\nu}}$ requires an increase in $\bar{\tau}_{\nu}$.

Finally, Figure~\ref{Fig:AugImp}(b) plots the system norm as a function of $\rho$ for the three augmented algorithms, for a test case with $n=4$ agents. The norm is not a monotonic function of the augmentation factor $\rho$ for the implementations in $\mathrm{RA}_{\rm cent}(\rho)$ and $\mathrm{RA}_{\rm dist}(\rho)$, but is monotonic for $\mathrm{RA}_{\rm dist}^{\rm dual}(\rho)$ applied to resource allocation problems, in agreement with the result for the parametrically uniform case in equation \eqref{Eq:Whatever}.

\section{Conclusions}
\label{Section: Conclusions}

We have studied the input-output performance of continuous-time saddle-point methods for solving linearly constrained convex quadratic programs, providing an explicit formula for the $\mathcal{H}_2$ norm under a relevant input-output configuration. We then studied the effects of Lagrangian regularization and augmentation on this norm, and derived a distributed dual version of the augmented algorithm which overcomes some of the limitations of naive augmentation. We then applied the results to compare several implementations of the saddle-point method to resource allocation problems.
Open directions for future research include input-output performance metrics for problems involving inequality constraints, for distributed implementations where communication delays occur between agents, and for other classes of distributed optimization algorithms. An analogous study in a $\mathcal{H}_{\infty}$ performance framework has not been completed. Another interesting question is how to further improve the $\mathcal{H}_2$ performance of saddle-point methods by designing auxiliary feedback controllers; augmentation is but one approach. Finally, extending these results beyond the case of quadratic cost functions with linear constraints will require nonlinear, robust control approaches, as in \cite{AC-EM-SL-JC:16,AC-EM-SL-JC:17,JWSP:16f,LL-BR-AP:16,BH-LL:17}.


\renewcommand{\baselinestretch}{1}
\bibliographystyle{IEEEtran}

\begin{IEEEbiography}[{\includegraphics[width=1in,height=1.25in,clip,keepaspectratio]{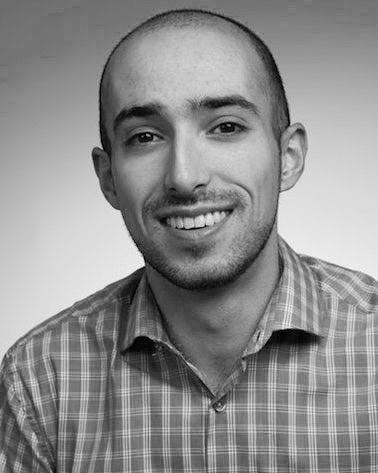}}]{John W. Simpson-Porco} (S'11--M'16) received the B.Sc. degree in engineering physics from Queen's University, Kingston, ON, Canada in 2010, and the Ph.D. degree in mechanical engineering from the University of California at Santa Barbara, Santa Barbara, CA, USA in 2015.

He is currently an Assistant Professor of Electrical and Computer Engineering at the University of Waterloo, Waterloo, ON, Canada. He was previously a visiting scientist with the Automatic Control Laboratory at ETH Z\"{u}rich, Z\"{u}rich, Switzerland. His research focuses on the control and optimization of multi-agent systems and networks, with applications in modernized power grids.

Prof. Simpson-Porco is a recipient of the 2012--2014 IFAC Automatica Prize and the Center for Control, Dynamical Systems and Computation Best Thesis Award and Outstanding Scholar Fellowship.
\end{IEEEbiography}

\begin{IEEEbiography}%
[{\includegraphics[width=1in,height=1.45in,clip,keepaspectratio]{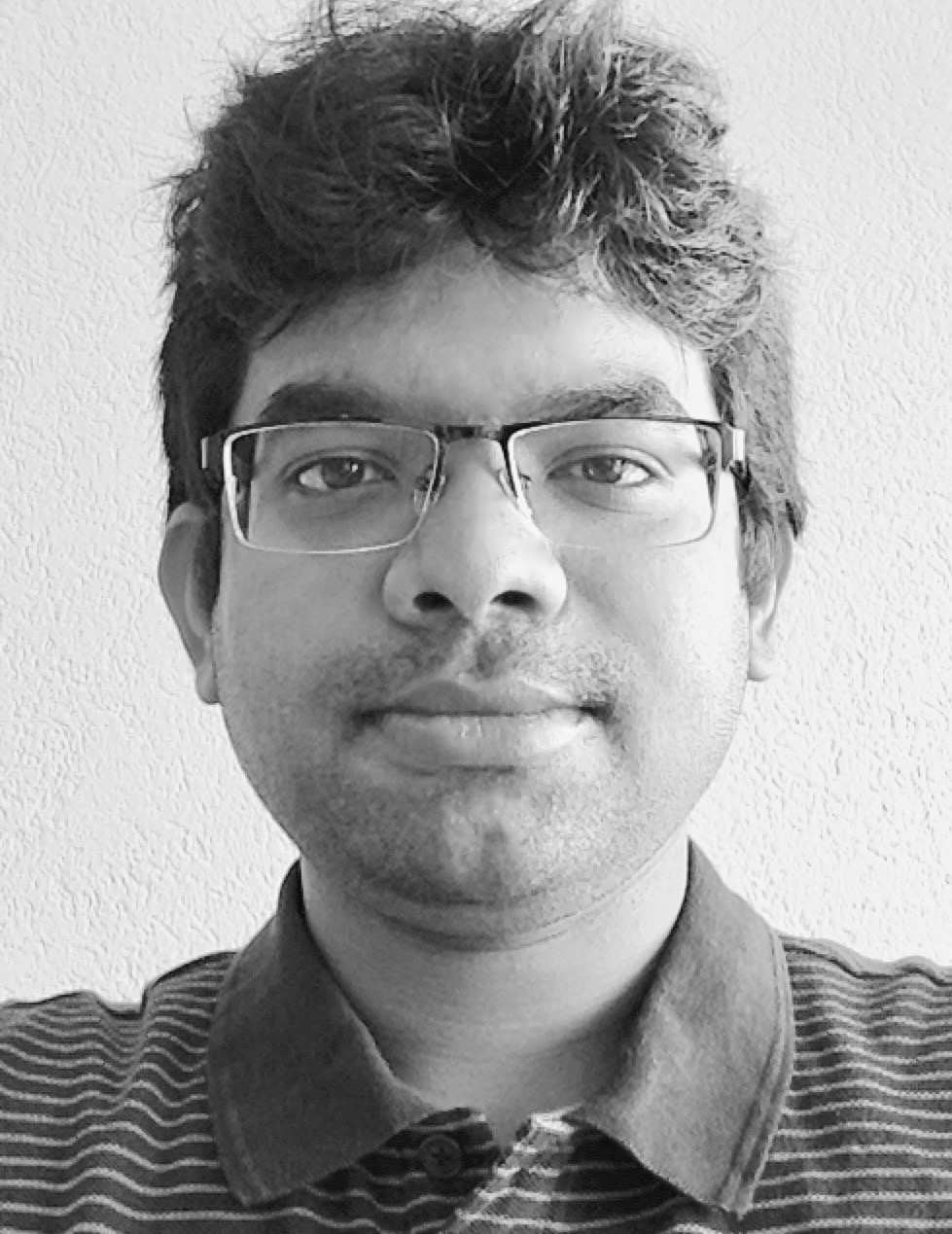}}]%
{Bala Kameshwar Poolla} (S'15) received his B.Tech. degree in Electrical Engineering and M.Tech. degree in Control Systems Engineering from the Indian Institute of Technology Kharagpur, India in 2013. He is currently a Ph.D. student at the Automatic Control Laboratory at ETH Z\"{u}rich. His research interests include applications of control theory to power grids, renewable energy systems, and energy markets. He was a finalist for the ACC 2016 Best Student Paper Award.
\end{IEEEbiography}

\begin{IEEEbiography}%
[{\includegraphics[width=1in,height=1.25in,clip,keepaspectratio]{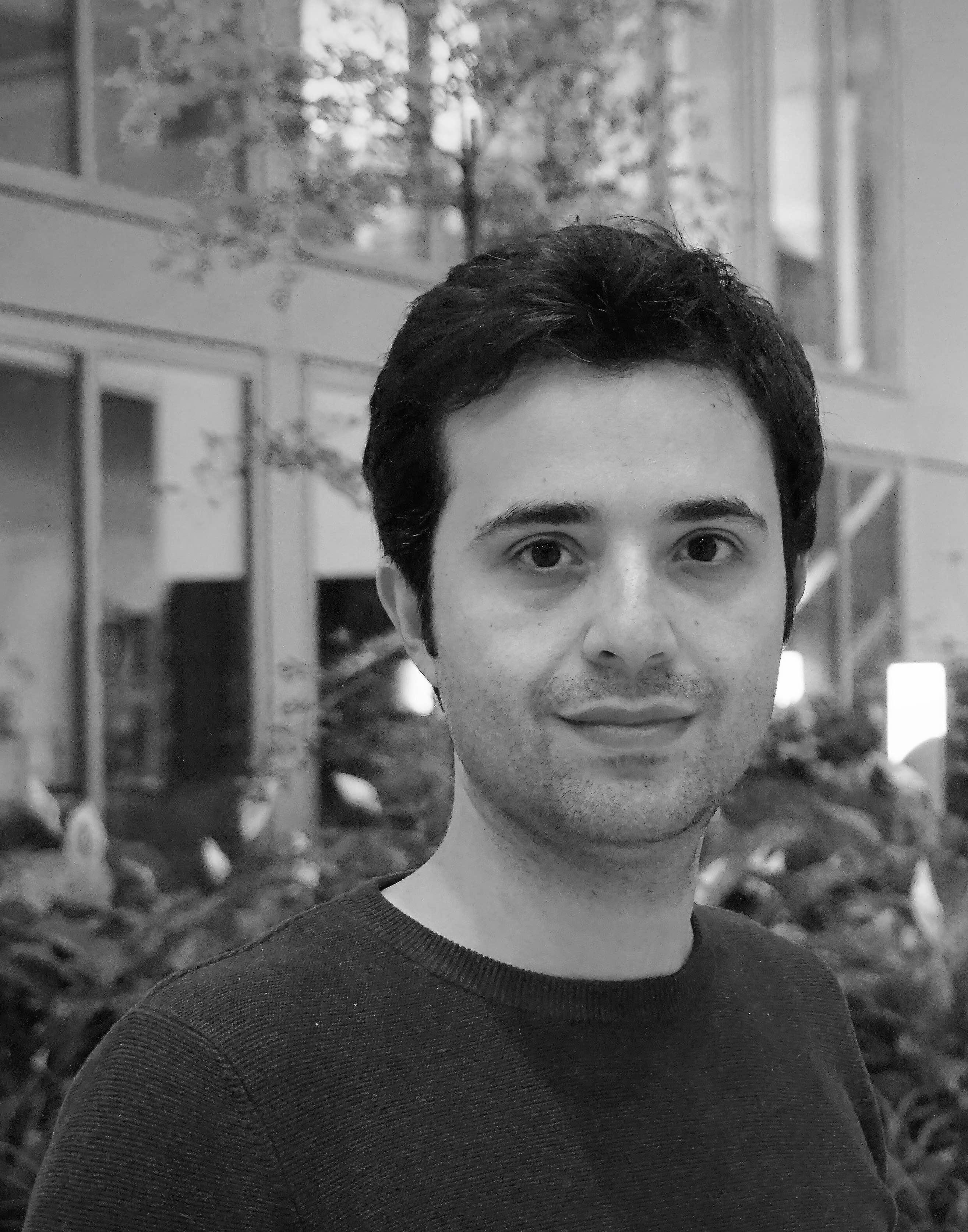}}]%
{Nima Monshizadeh} obtained his B.Sc. degree in electrical engineering from the University of Tehran, the M.Sc. degree in control engineering from K.N. Toosi University of Technology, Tehran, Iran, and the Ph.D. degree with honor (cum laude) from the Johann Bernoulli Institute for Mathematics and Computer Science, University of Groningen, The Netherlands, in December 2013. Afterwards, he was appointed as a Post-Doctoral Researcher at the Engineering and Technology Institute of the University of Groningen. He was a Research Associate with the Control group of the University of Cambridge, before taking a tenure track position at the University of Groningen in January 2018. His research interests include power networks, model reduction, optimization and control of complex networks.
\end{IEEEbiography}

\begin{IEEEbiography}%
[{\includegraphics[width=1in,height=1.25in,clip,keepaspectratio]{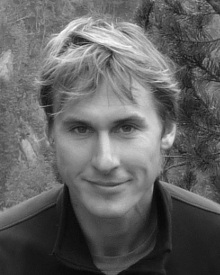}}]%
{Florian D\"{o}rfler} (S'09--M'13) is an Assistant Professor at the Automatic Control Laboratory at ETH Z\"urich. He received his Ph.D. degree in Mechanical Engineering from the University of California at Santa Barbara in 2013, and a Diplom degree in Engineering Cybernetics from the University of Stuttgart in 2008. From 2013 to 2014 he was an Assistant Professor at the University of California Los Angeles. His primary research interests are centered around distributed control, complex networks, and cyber-physical systems currently with  applications in energy systems and smart grids. His students were winners or finalists for Best Student Paper awards at the European Control Conference (2013), the American Control Conference (2016), and the PES PowerTech Conference 2017. His articles received the 2010 ACC Student Best Paper Award, the 2011 O. Hugo Schuck Best Paper Award, and the 2012-2014 Automatica Best Paper Award, and the 2016 IEEE Circuits and Systems Guillemin-Cauer Best Paper Award. He is a recipient of the 2009 Regents Special International Fellowship, the 2011 Peter J. Frenkel Foundation Fellowship, and the 2015 UCSB ME Best PhD award.

\end{IEEEbiography}

\end{document}